# Differential equation approximations for Markov chains

**R.W.R. Darling**

*Mathematics Research Group, National Security Agency, 9800 Savage Road
Fort George G. Meade, Maryland 20755-6515, USA
e-mail:* rwdarli@nsa.gov

and

**J.R. Norris**

*University of Cambridge, Statistical Laboratory, Centre for Mathematical Sciences
Wilberforce Road, Cambridge, CB3 0WB, UK
e-mail:* j.r.norris@statslab.cam.ac.uk

**Abstract:** We formulate some simple conditions under which a Markov chain may be approximated by the solution to a differential equation, with quantifiable error probabilities. The role of a choice of coordinate functions for the Markov chain is emphasised. The general theory is illustrated in three examples: the classical stochastic epidemic, a population process model with fast and slow variables, and core-finding algorithms for large random hypergraphs.

## Contents

## 1. Introduction

Differential equations and Markov chains are the basic models of dynamical systems in a deterministic and a probabilistic context, respectively. Since the analysis of differential equations is often more feasible and efficient, both from a mathematical and a computational point of view, it is of interest to understand in some generality when the sample paths of a Markov chain can be guaranteed to lie, with high probability, close to the solution of a differential equation.

We shall obtain a number of estimates, given explicitly in terms of the Markov transition rates, for the probability that a Markov chain deviates further than a given distance from the solution to a suitably chosen differential equation. The basic method is simply a combination of Gronwall's lemma with martingale inequalities. The intended contribution of this paper is to set out in a convenient form some estimates that can be deduced in this way, along with some illustrations of their use. Although it is widely understood how to arrive at a suitable differential equation, the justification of an approximation statement can be more challenging, particularly if one has cause to push beyond the scope of classical weak convergence results. We have found the use of explicit estimates effective, for example, when the Markov chain terminates abruptly on leaving some domain [4], or when convergence is needed over a long time interval [23], or for processes having a large number of components with very different scales [14].

The first step in our approach is a choice of coordinate functions for the given Markov chain: these are used to rescale the process, whose values might typically form a vector of non-negative integers, to one which may lie close to a continuously evolving path. The choice of coordinate functions may also be used to forget some components of the Markov chain which do not behave suitably and further, as is sometimes necessary, to correct the values of the remaining components to take account of the values of the forgotten components. This is illustrated in the examples in Sections 6 and 7. The behaviour of forgotten components can sometimes be approximated by a random process having relatively simple characteristics, which are determined by the differential equation. This is illustrated in the example in Section 5, where it is used to show the asymptotic independence of individuals in a large population.

We have been motivated by two main areas of application. The first is to population processes, encompassing epidemic models, queueing and network models, and models for chemical reactions. It is often found in models of interest that certain variables oscillate rapidly and randomly while others, suitably rescaled, are close to deterministic. It was a primary motivation to find an extension of our quantitative estimates which was useful in such a context. The example in Section 6, which is drawn from [2], shows that this is possible. The second area of application is the analysis of randomized algorithms and combinatorial structures. Here, the use of differential equation approximations has become an important tool. The example in Section 7 gives an alternative treatment and generalization of the $k$-core asymptotics discovered in [19].

The martingale estimates we need are derived from scratch in the Appendix, using a general procedure for the identification of martingales associated to a Markov chain. We have taken the opportunity to give a justification of this procedure, starting from a presentation of the chain in terms of its jump chain and holding times. We found it interesting to do this without passing through the characterization of Markov chains in terms of semigroups and generators.

The authors are grateful to Perla Sousi and to a referee for a careful reading of an earlier version of this paper, which has helped to clarify the present work.

## 2. Survey of related literature

There is a well-developed body of literature devoted to the general question of the convergence of Markov processes, which includes as a special case the question we address in this paper. This special case arises under *fluid limit* or *law of large numbers scaling*, where, for large $N$, jumps whose size is of order $1/N$ occur at a rate of order $N$. This is to be distinguished from *diffusive* or *central limit scaling*, where jumps of mean zero and of size of order $1/\sqrt{N}$ occur at a rate of order $N$. Just as in the classical central limit theorem, a Gaussian diffusive limit can be used to describe to first order the fluctuations of a process around its fluid limit.

Both sorts of limit are presented in the books by Ethier and Kurtz [6, Section 7.4], Jacod and Shiryaev [8, Section IX.4b], and Kallenberg [9]. These works develop conditions on the transition operators or rate kernels of a sequence of Markov chains which are sufficient to imply the weak convergence of the corresponding processes. Trotter's paper [22] was one of the first to take this point of view.

The fluid limit is more elementary and often allows, with advantage, a more direct approach. One identifies a limiting drift $b$ of the processes, which we shall suppose to be a Lipschitz vector field, and then the limit is the deterministic path obtained by solving the differential equation $\dot x = b(x)$. Kurtz [12] describes some sufficient conditions for weak convergence of processes in this context. Since the limit is continuous in time, weak convergence is here simply convergence in probability to 0 of the maximal deviation from the limit path over any given compact time interval. Later, exponential martingale estimates, were used to prove decay of error probabilities at an exponential rate. See the book of Shwartz and Weiss [21]. This is the direction also of the present paper. Differential equation approximations for stochastic systems with small noise have been studied for many sorts of process other than Markov processes. See the book of Kushner and Yin [13].

Applications of fluid limits for Markov chains are scattered across many fields. See [3] on epidemiology and [21] on communications and computer networks. Much has been achieved by the identification of deterministic limit behaviour when randomized algorithms are applied to large combinatorial problems, or deterministic algorithms are applied to large random combinatorial structures. Examples include Karp and Sipser's seminal paper [10] on maximal matchings, Hajek's analysis [7] of communications protocols, Mitzenmacher's [16] balanced allocations, and the analysis of Boolean satisfiability by Achlioptas [1] and Semerjian and Monasson [18]. A general framework for this sort of application was developed by Wormald and others, see [19], [24].

Finally, the emergence of deterministic macroscopic evolutions from microsopic behaviour, often assumed stochastic, is a more general phenomenon than addressed in the literature mentioned above. We have only considered scaling the sizes of the components of a Markov chain. In random models where each component counts the number of particles at a given spatial location, it is natural to scale also these locations, leading sometimes to macroscopic laws governed by partial rather than ordinary differential equations. This is the field of hydrodynamic limits – see, for example, Kipnis and Landim [11], for an introduction.

## 3. Some simple motivating examples

We now give a series of examples of Markov processes, each of which takes many small jumps at a fast rate. The *drift* is the product of the average jump by the total rate, which may vary from state to state. In cases where there are a number of different types of jump, one can compute the drift as a sum over types of the size of the jump multiplied by its rate. We write down the drift and hence obtain a differential equation. In the rest of the paper, we give conditions under which the Markov chain will be well approximated by solutions of this equation. In each of the examples there is a parameter $N$ which quantifies the smallness of the jumps and the compensating largeness of the jump rates. The approximations will be good when $N$ is large.

### 3.1. Poisson process

Take $(X_t)_{t\geqslant 0}$ to be a Poisson process of rate $\lambda N$, and set $\boldsymbol{X}_t = X_t/N$. Note that $\boldsymbol{X}$ takes jumps of size $1/N$ at rate $\lambda N$. The drift is then $\lambda$ and the differential equation is
$$\dot{x}_t = \lambda.$$
If we take as initial state $\boldsymbol{X}_0 = x_0 = 0$, then we may expect that $\boldsymbol{X}_t$ stay close to the solution $x_t = \lambda t$. This is a law of large numbers for the Poisson process.

### 3.2. $M_N/M_1/\infty$ queue

Consider a queue with arrivals at rate $N$, exponential service times of mean 1, and infinitely many servers. Write $X_t$ for the number of customers present at time $t$. Set $\boldsymbol{X}_t = X_t/N$, then $\boldsymbol{X}$ is a Markov chain, which jumps by $1/N$ at rate $N$, and jumps by $-1/N$ at rate $N\boldsymbol{X}_t$. The drift is then $1-x$ and the differential equation is
$$\dot{x}_t = 1 - x_t.$$
The solution of this equation is given by $x_t = 1 + x_0 e^{-t}$, so we may expect that, for large $N$ the queue size stabilizes near $N$, at exponential rate 1.

### 3.3. Chemical reaction $A + B \leftrightarrow C$

In a reversible reaction, pairs of molecules of types $A$ and $B$ become a single molecule of type $C$ at rate $\lambda/N$, and molecules of type $C$ become a pair of molecules of types $A$ and $B$ at rate $\mu$. Write $A_t, B_t, C_t$ for the numbers of molecules of each type present at time $t$. Set
$$\boldsymbol{X}_t = (X_t^1, X_t^2, X_t^3) = (A_t, B_t, C_t)/N,$$
then $\boldsymbol{X}$ is a Markov chain, which makes jumps of $(-1,-1,1)/N$ at rate $(\lambda/N)(NX_t^1)(NX_t^2)$, and makes jumps of $(1,1,-1)/N$ at rate $\mu(NX_t^3)$. The drift is then $(\mu x^3 - \lambda x^1 x^2, \mu x^3 - \lambda x^1 x^2, \lambda x^1 x^2 - \mu x^3)$ and the differential equation is, in components,
$$\dot{x}_t^1 = \mu x_t^3 - \lambda x_t^1 x_t^2, \quad \dot{x}_t^2 = \mu x_t^3 - \lambda x_t^1 x_t^2, \quad \dot{x}_t^3 = \lambda x_t^1 x_t^2 - \mu x_t^3.$$
Any vector $(x^1, x^2, x^3)$ with $\mu x^3 = \lambda x^1 x^2$ is a fixed point of this equation and may be expected to correspond to an equilibrium state of the system.

### 3.4. Gunfight

Two gangs of gunmen fire at each other. On each side, each surviving gunman hits one of the opposing gang randomly, at rate $\alpha$ for gang $A$ and at rate $\beta$ for gang $B$. Write $A_t$ and $B_t$ for the numbers still firing on each side at time $t$. Set $\boldsymbol{X}_t = (X_t^1, X_t^2) = (A_t, B_t)/N$, then $\boldsymbol{X}$ is a Markov chain and jumps by $(0, -1)/N$ at rate $\alpha N X_t^1$, and by $(-1, 0)/N$ at rate $\beta N X_t^2$. The drift is then $(-\beta x^2, -\alpha x^1)$ and the differential equation is, in components,

$$\dot{x}_t^1 = -\alpha x_t^2, \quad \dot{x}_t^2 = -\beta x_t^1.$$

Note that, in this case the parameter $N$ does not enter the description of the model. However the theory will give an informative approximation only for initial conditions of the type $(A_t, B_t) = N(a_0, b_0)$. The reader may like to solve the equation and see who wins the fight.

### 3.5. Continuous time branching processes

Each individual in a population lives for an exponentially distributed time of mean $1/N$, whereupon it is replaced by a random number $Z$ of identical offspring, where $Z$ has finite mean $\mu$. Distinct individuals behave independently. Write $X_t$ for the number of individuals present at time $t$. Set $\boldsymbol{X}_t = X_t/N$, then $\boldsymbol{X}$ is a Markov chain, which jumps by $(k-1)/N$ at rate $N\boldsymbol{X}_t \mathbb{P}(Z = k)$ for all $k \in \mathbb{Z}^+$. The drift is then $\sum_k (k-1)\mathbb{P}(Z = k)x = (\mu - 1)x$ and the differential equation is

$$\dot{x}_t = (\mu - 1)x_t.$$

This equation gives a first order approximation for the evolution of the population size – in particular, it is clear that the cases where $\mu < 1$, $\mu = 1$, $\mu > 1$ should show very different long-time behaviour.

## 4. Derivation of the estimates

Let $X = (X_t)_{t \geqslant 0}$ be a continuous-time Markov chain with countable[1] state-space $S$. Assume that in every state $\xi \in S$ the total jump rate $q(\xi)$ is finite, and write $q(\xi, \xi')$ for the jump rate from $\xi$ to $\xi'$, for each pair of distinct states $\xi$ and $\xi'$. We assume that $X$ does not explode: a simple sufficient condition for this is that the jump rates are bounded, another is that $X$ is recurrent.

We make a choice of *coordinate functions* $x^i : S \to \mathbb{R}$, for $i = 1, \ldots, d$, and write $\boldsymbol{x} = (x^1, \ldots, x^d) : S \to \mathbb{R}^d$. Consider the $\mathbb{R}^d$-valued process $\boldsymbol{X} = (\boldsymbol{X}_t)_{t \geqslant 0}$ given by $\boldsymbol{X}_t = (X_t^1, \ldots, X_t^d) = \boldsymbol{x}(X_t)$. Define, for each $\xi \in S$, the *drift vector*

$$\beta(\xi) = \sum_{\xi' \neq \xi} (\boldsymbol{x}(\xi') - \boldsymbol{x}(\xi))q(\xi, \xi'),$$

where we set $\beta(\xi) = \infty$ if this sum fails to converge absolutely.

---
[1] The extension of the results of this paper to the case of a general measurable state-space is a routine exercise.

Our main goal is the derivation of explicit estimates which may allow the approximation of $\boldsymbol{X}$ by the solution of a differential equation. We shall also discuss how the computation of certain associated probabilities can be simplified when such an approximation is possible.

Let $U$ be a subset of $\mathbb{R}^d$ and let $x_0 \in U$. Let $b : U \to \mathbb{R}^d$ be a Lipschitz vector field. The differential equation $\dot{x}_t = b(x_t)$ has a unique maximal solution $(x_t)_{t \leqslant \zeta}$, starting from $x_0$, with $x_t \in U$ for all $t < \zeta$. Maximal here refers to $\zeta$ and means that there is no solution in $U$ defined on a longer time interval. Our analysis is based on a comparison of the equations

$$\boldsymbol{X}_t = \boldsymbol{X}_0 + M_t + \int_0^t \beta(X_s)ds, \quad 0 \leqslant t \leqslant T_1,$$

$$x_t = x_0 + \int_0^t b(x_s)ds, \quad 0 \leqslant t \leqslant \zeta,$$

where $T_1 = \inf\{t \geqslant 0 : \beta(X_t) = \infty\}$ and where the first equation serves to define the process $(M_t)_{0 \leqslant t \leqslant T_1}$.

## 4.1. $L^2$-estimates

The simplest estimate we shall give is obtained by a combination of Doob's $L^2$-inequality and Gronwall's lemma. Doob's $L^2$-inequality states that, for any martingale $(M_t)_{t \leqslant t_0}$,

$$\mathbb{E}\left(\sup_{t \leqslant t_0} |M_t|^2\right) \leqslant 4\mathbb{E}\left(|M_{t_0}|^2\right).$$

Gronwall's lemma states that, for any real-valued integrable function $f$ on the interval $[0, t_0]$, the inequality

$$f(t) \leqslant C + D \int_0^t f(s)ds, \quad \text{for all } t, \tag{1}$$

implies that $f(t_0) \leqslant Ce^{Dt_0}$.

Write, for now, $K$ for the Lipschitz constant of $b$ on $U$ with respect to the Euclidean norm $|.|$. Fix $t_0 < \zeta$ and $\varepsilon > 0$ and assume that[2],

$$\text{for all } \xi \in S \text{ and } t \leqslant t_0, \quad |\boldsymbol{x}(\xi) - x_t| \leqslant \varepsilon \implies \boldsymbol{x}(\xi) \in U.$$

Set $\delta = \varepsilon e^{-Kt_0}/3$ and fix $A > 0$. For our estimate to be useful it will be necessary that $A$ be small compared to $\varepsilon^2$. Set $T = \inf\{t \geqslant 0 : \boldsymbol{X}_t \notin U\}$; see Figure 1. Define, for $\xi \in S$,

$$\alpha(\xi) = \sum_{\xi' \neq \xi} |\boldsymbol{x}(\xi') - \boldsymbol{x}(\xi)|^2 q(\xi, \xi').$$

---
[2] A simpler but stronger condition is to require that path $(x_t)_{t \leqslant t_0}$ lies at a distance greater than $\varepsilon$ from the complement of $U$.

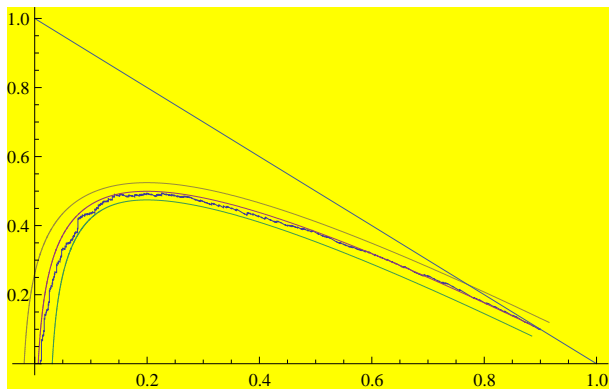

FIG 1. *The unit square $(0, 1)^2$ is a possible choice of the set $U$. The inner and outer solid curves bound a tube around the deterministic solution (red) of the ordinary differential equation, which starts inside $U$. The realization shown of the Markov chain trajectory does not leave the tube before exit from $U$. This is a realization of the stochastic epidemic, which will be discussed in more detail in Section 5.*

Consider the events[3]

$$\Omega_0 = \{|\boldsymbol{X}_0 - x_0| \leqslant \delta\}, \quad \Omega_1 = \left\{\int_0^{T \wedge t_0} |\beta(X_t) - b(\boldsymbol{x}(X_t))| dt \leqslant \delta\right\}$$

and

$$\Omega_2 = \left\{\int_0^{T \wedge t_0} \alpha(X_t) dt \leqslant At_0\right\}, \quad \Omega_2' = \left\{T \wedge t_0 \leqslant T_1 \text{ and } \sup_{t \leqslant T \wedge t_0} |M_t| \leqslant \delta\right\}.$$

Consider the random function $f(t) = \sup_{s \leqslant t} |\boldsymbol{X}_s - x_s|$ on the interval $[0, T \wedge t_0]$. Then

$$f(t) \leqslant |\boldsymbol{X}_0 - x_0| + \sup_{s \leqslant t} |M_s| + \int_0^t |\beta(X_s) - b(\boldsymbol{x}(X_s))| ds + \int_0^t |b(\boldsymbol{X}_s) - b(x_s)| ds.$$

So, on the event $\Omega_0 \cap \Omega_1 \cap \Omega_2'$, $f$ satisfies (1) with $C = 3\delta$ and $D = K$, so $f(T \wedge t_0) \leqslant \varepsilon$, which implies $T > t_0$ and hence $f(t_0) \leqslant \varepsilon$. Consider now the stopping time

$$\tilde{T} = T \wedge t_0 \wedge \inf\left\{t \geqslant 0 : \int_0^t \alpha(X_s) ds > At_0\right\}.$$

---

[3]In examples, we shall often have some or all of these events equal to $\Omega$. We may have $\boldsymbol{X}_0 = x_0$ and $\beta = b \circ \boldsymbol{x}$, or at least be able to show that $\beta - b \circ \boldsymbol{x}$ is uniformly small on $S$. The example discussed in Subsection 6.1 exploits fully the given form of $\Omega_1$, in that the integrand $\beta(X_t) - b(\boldsymbol{x}(X_t))$ cannot be bounded uniformly in a suitable way, whereas the integral is suitably small with high probability.

By Cauchy–Schwarz, we have $|\beta(\xi)|^2 \leqslant q(\xi)\alpha(\xi)$ for all $\xi \in S$, so $\tilde{T} \leqslant T_1$. By a standard argument using Doob's $L^2$-inequality, which is recalled in Proposition 8.7, we have

$$\mathbb{E}\left(\sup_{t \leqslant \tilde{T}} |M_t|^2\right) \leqslant 4At_0.$$

On $\Omega_2$, we have $\tilde{T} = T \wedge t_0$, so $\Omega_2 \setminus \Omega_2' \subseteq \{\sup_{t \leqslant \tilde{T}} |M_t| > \delta\}$ and so, by Chebyshev's inequality, $\mathbb{P}(\Omega_2 \setminus \Omega_2') \leqslant 4At_0/\delta^2$. We have proved the following result, which can sometimes enable us to show that the situation illustrated in Figure 1 occurs with high probability.

**Theorem 4.1.** *Under the above conditions,*

$$\mathbb{P}\left(\sup_{t \leqslant t_0} |\boldsymbol{X}_t - x_t| > \varepsilon\right) \leqslant 4At_0/\delta^2 + \mathbb{P}(\Omega_0^c \cup \Omega_1^c \cup \Omega_2^c).$$

## *4.2. Exponential estimates*

It is clear that the preceding argument could be applied for any norm on $\mathbb{R}^d$ with obvious modifications. We shall do this for the supremum norm $\|x\| = \max_i |x_i|$, making at the same time a second variation in replacing the use of Doob's $L^2$-inequality with an exponential martingale inequality. This leads to the version of the result which we prefer for the applications we have considered. It will be necessary to modify some assumptions and notation introduced in the preceding subsection. We shall stick to these modifications from now on. We assume now that, $\varepsilon > 0$ and $t_0$ are chosen so that,

$$\text{for all } \xi \in S \text{ and } t \leqslant t_0, \quad \|\boldsymbol{x}(\xi) - x_t\| \leqslant \varepsilon \implies \boldsymbol{x}(\xi) \in U. \tag{2}$$

Write now $K$ for the Lipschitz constant of $b$ with respect to the supremum norm. Set $\delta = \varepsilon e^{-Kt_0}/3$. Fix $A > 0$ and set $\theta = \delta/(At_0)$. Define

$$\sigma_\theta(x) = e^{\theta|x|} - 1 - \theta|x|, \quad x \in \mathbb{R}$$

and set

$$\phi^i(\xi, \theta) = \sum_{\xi' \neq \xi} \sigma_\theta(x^i(\xi') - x^i(\xi))q(\xi, \xi'), \quad \phi(\xi, \theta) = \max_i \phi^i(\xi, \theta), \quad \xi \in S.$$

Consider the events

$$\Omega_0 = \{\|\boldsymbol{X}_0 - x_0\| \leqslant \delta\}, \quad \Omega_1 = \left\{\int_0^{T \wedge t_0} \|\beta(X_t) - b(\boldsymbol{x}(X_t))\| dt \leqslant \delta\right\},$$

$$\Omega_2 = \left\{\int_0^{T \wedge t_0} \phi(X_t, \theta) dt \leqslant \tfrac{1}{2}\theta^2 At_0\right\}$$

and
$$\Omega_2' = \left\{ T \wedge t_0 \leqslant T_1 \text{ and } \sup_{t \leqslant T \wedge t_0} \|M_t\| \leqslant \delta \right\}.$$

We can use Gronwall's lemma, as above, to see that on the event $\Omega_0 \cap \Omega_1 \cap \Omega_2'$ we have $\sup_{t \leqslant t_0} \|\boldsymbol{X}_t - x_t\| \leqslant \varepsilon$. Note that, since $\sigma_\theta(x) \geqslant \theta^2 |x|^2/2$ for all $x \in \mathbb{R}$, we always have $T \wedge t_0 \leqslant T_1$ on $\Omega_2$. Fix $i \in \{1, \ldots, d\}$ and set

$$\phi(\xi) = \sum_{\xi' \neq \xi} \left\{ e^{\theta(x^i(\xi') - x^i(\xi))} - 1 - \theta(x^i(\xi') - x^i(\xi)) \right\}.$$

Then $\phi(\xi) \leqslant \phi^i(\xi, \theta) \leqslant \phi(\xi, \theta)$, so

$$\mathbb{P}\left( \sup_{t \leqslant T \wedge t_0} M_t^i > \delta \text{ and } \Omega_2 \right) \leqslant \mathbb{P}\left( \sup_{t \leqslant T \wedge t_0} M_t^i > \delta \text{ and } \int_0^{T \wedge t_0} \phi(X_t) dt \leqslant \tfrac{1}{2}\theta^2 A t_0 \right)$$
$$\leqslant \exp\{\tfrac{1}{2}\theta^2 A t_0 - \theta\delta\} = \exp\{-\delta^2/(2At_0)\}.$$

For the second inequality, we used a standard exponential martingale inequality, which is recalled in Proposition 8.8. Since the same argument applies also to $-M$ and for all $i$, we thus obtain $\mathbb{P}(\Omega_2 \setminus \Omega_2') \leqslant 2de^{-\delta^2/(2At_0)}$. We have proved the following estimate, which is often stronger than Theorem 4.1. In an asymptotic regime where the sizes of jumps in $\boldsymbol{X}$ are of order $1/N$ but their rates are of order $N$, the estimate will often allow us to prove decay of error probabilities in the differential equation approximation at a rate exponential in $N$. The price to be paid for this improvement is the necessity to deal with the event $\Omega_2$ just defined rather than its more straightforward counterpart in the preceding subsection[4].

**Theorem 4.2.** *Under the above conditions,*
$$\mathbb{P}\left( \sup_{t \leqslant t_0} \|\boldsymbol{X}_t - x_t\| > \varepsilon \right) \leqslant 2de^{-\delta^2/(2At_0)} + \mathbb{P}(\Omega_0^c \cup \Omega_1^c \cup \Omega_2^c).$$

### 4.3. Convergence of terminal values

In cases where the solution of the differential equation leaves $U$ in a finite time, so that $\zeta < \infty$, we can adapt the argument to obtain estimates on the time $T$ that $\boldsymbol{X}$ leaves $U$ and on the terminal value $\boldsymbol{X}_T$. The vector field $b$ can be extended to the whole of $\mathbb{R}^d$ with the same Lipschitz constant. Let us choose such an extension, also denoted $b$, and write now $(x_t)_{t \geqslant 0}$ for the unique solution to $\dot{x}_t = b(x_t)$ starting from $x_0$. Define for $\varepsilon > 0$

$$\zeta_\varepsilon^- = \inf\{t \geqslant 0 : x \notin U \text{ for some } x \in \mathbb{R}^d \text{ with } \|x - x_t\| \leqslant \varepsilon\},$$
$$\zeta_\varepsilon^+ = \inf\{t \geqslant 0 : x \notin U \text{ for all } x \in \mathbb{R}^d \text{ with } \|x - x_t\| \leqslant \varepsilon\}.$$

---

[4]The present approach is useful only when the jumps of $\boldsymbol{X}$ have an exponential moment, whereas the previous $L^2$ approach required only jumps of finite variance. In many applications, the jumps are uniformly bounded: if $J$ is an upper bound for the supremum norm of the jumps, then, using the inequality $e^x - 1 - x \leqslant \tfrac{1}{2}x^2 e^x$, a sufficient condition for $\Omega_2 = \Omega$ is that $A \geqslant QJ^2 \exp\{\delta J/(At_0)\}$, where $Q$ is the maximum jump rate.

and set[5]

$$\rho(\varepsilon) = \sup\{\|x_t - x_\zeta\| : \zeta_\varepsilon^- \leqslant t \leqslant \zeta_\varepsilon^+\}.$$

Typically we will have $\rho(\varepsilon) \to 0$ as $\varepsilon \to 0$. Indeed, if $U$ has a smooth boundary at $x_\zeta$, and if $b(x_\zeta)$ is not tangent to this boundary, then $\rho(\varepsilon) \leqslant C\varepsilon$ for some $C < \infty$, for all sufficiently small $\varepsilon > 0$. However, we leave this step until we consider specific examples. Assume now, in place of (2), that $\varepsilon$ and $t_0$ are chosen so that $t_0 > \zeta_\varepsilon^+$. On $\Omega_0 \cap \Omega_1 \cap \Omega_2'$, we obtain, as above, that $f(T \wedge t_0) \leqslant \varepsilon$, which forces $\zeta_\varepsilon^- \leqslant T \leqslant \zeta_\varepsilon^+$, and hence $\|X_T - x_\zeta\| \leqslant \|X_T - x_T\| + \|x_T - x_\zeta\| \leqslant \varepsilon + \rho(\varepsilon)$. We have proved the following estimate[6].

**Theorem 4.3.** *Under the above conditions,*

$$\mathbb{P}\big(\|X_T - x_\zeta\| > \varepsilon + \rho(\varepsilon)\big) \leqslant \mathbb{P}\left(T \notin [\zeta_\varepsilon^-, \zeta_\varepsilon^+] \text{ or } \sup_{t \leqslant T}\|X_t - x_t\| > \varepsilon\right)$$
$$\leqslant 2de^{-\delta^2/(2At_0)} + \mathbb{P}(\Omega_0^c \cup \Omega_1^c \cup \Omega_2^c).$$

### 4.4. Random processes modulated by the fluid limit

We return now to the case where $t_0 < \zeta$ and condition (2) holds. Although the results given so far can be interpreted as saying that $X$ is close to deterministic, there are sometimes associated random quantities which we may wish to understand, and whose behaviour can be described, approximately, in a relatively simple way in terms of the deterministic path $(x_t)_{t \leqslant t_0}$. To consider this in some generality, suppose there is given a countable set $I$ and a function $y : S \to I$ and consider the process $Y = (Y_t)_{t \geqslant 0}$ given by $Y_t = y(X_t)$. Define, for $\xi \in S$ and $y \in I$ with $y \neq y(\xi)$, the jump rates

$$\gamma(\xi, y) = \sum_{\xi' \in S : y(\xi') = y} q(\xi, \xi').$$

We now give conditions which may allow us to approximate $Y$ by a Markov chain with time-dependent jump rates, which are given in terms of the path $(x_t)_{t \leqslant t_0}$ and a non-negative function $g$ on $U \times \{(y, y') \in I \times I : y \neq y'\}$. Set $g_t(y, y') = g(x_t, y, y')$ for $t \leqslant t_0$. Fix $I_0 \subseteq I$ and set

$$\kappa = \sup_{t \leqslant t_0} \sup_{\|x - x_t\| \leqslant \varepsilon, y \in I_0} \sum_{y' \neq y} |g(x, y, y') - g(x_t, y, y')|.$$

---

[5] The function $\rho$ depends on the choice of extension made of $b$ outside $U$, whereas the distribution of $\|X_T - x_\zeta\|$ does not. This is untidy, but it is not simple to optimise over Lipschitz extensions, and in any case, this undesirable dependence of $\rho$ is a second order effect as $\varepsilon \to 0$.

[6] The same argument can be made using, in place of $\zeta_\varepsilon^\pm$, the times

$$\tilde{\zeta}_\varepsilon^- = \inf\{t \geqslant 0 : \boldsymbol{x}(\xi) \notin U \text{ for some } \xi \in S \text{ with } \|\boldsymbol{x}(\xi) - x_t\| \leqslant \varepsilon\},$$
$$\tilde{\zeta}_\varepsilon^+ = \inf\{t \geqslant 0 : \boldsymbol{x}(\xi) \notin U \text{ for all } \xi \in S \text{ with } \|\boldsymbol{x}(\xi) - x_t\| \leqslant \varepsilon\}.$$

This refinement can be useful if we wish to start $X$ on the boundary of $U$.

Set $T_0 = \inf\{t \geqslant 0 : X_t \notin U \text{ or } Y_t \notin I_0\}$, fix $G > 0$, and define

$$\Omega_3 = \left\{ \int_0^{T_0 \wedge t_0} \sum_{y \neq y(X_t)} |\gamma(X_t, y) - g(\boldsymbol{x}(X_t), y(X_t), y)| dt \leqslant G t_0 \right\}.$$

**Theorem 4.4.** *There exists a time-inhomogeneous Markov chain $(y_t)_{t \leqslant t_0}$, with state-space $I$ and jump rates $g_t(y, y')$, such that*

$$\mathbb{P}\left(\sup_{t \leqslant t_0} \|\boldsymbol{X}_t - x_t\| > \varepsilon \quad \text{or} \quad Y_t \neq y_t \text{ for some } t \leqslant \tau\right)$$

$$\leqslant (G + \kappa) t_0 + 2d e^{-\delta^2/(2At_0)} + \mathbb{P}(\Omega_0^c \cup \Omega_1^c \cup \Omega_2^c \cup \Omega_3^c),$$

*where $\tau = \inf\{t \geqslant 0 : y_t \notin I_0\} \wedge t_0$.*

*Proof.* We construct the process $(X_t, y_t)_{t \leqslant t_0}$ as a Markov chain, where the rates are chosen to keep the processes $(Y_t)_{t \leqslant t_0}$ and $(y_t)_{t \leqslant t_0}$ together for as long as possible. Define for $t \leqslant t_0$, and for $\xi, \xi' \in S$ and $y, y' \in I$, with $(\xi, y) \neq (\xi', y')$, first in the case $y = y(\xi)$,

$$q_t(\xi, y; \xi', y') = \begin{cases} q(\xi, \xi'), & \text{if } y' = y(\xi') = y(\xi), \\ q(\xi, \xi')\{1 \wedge (g_t(y, y')/\gamma(\xi, y'))\}, & \text{if } y' = y(\xi') \neq y(\xi), \\ q(\xi, \xi')\{1 - (g_t(y, y(\xi'))/\gamma(\xi, y(\xi')))\}^+, & \text{if } y' = y \neq y(\xi'), \\ \{g_t(y, y') - \gamma(\xi, y')\}^+, & \text{if } \xi' = \xi, \\ 0, & \text{otherwise,} \end{cases}$$

then in the case $y \neq y(\xi)$,

$$q_t(\xi, y; \xi', y') = \begin{cases} q(\xi, \xi'), & \text{if } y' = y, \\ g_t(y, y'), & \text{if } \xi' = \xi, \\ 0, & \text{otherwise.} \end{cases}$$

Consider the Markov chain $(X_t, y_t)_{t \leqslant t_0}$ on $S \times I$, starting from $(X_0, y(X_0))$, with jump rates $q_t(\xi, y; \xi', y')$. It is straightforward to check, by calculation of the marginal jump rates, that the components $(X_t)_{t \leqslant t_0}$ and $(y_t)_{t \leqslant t_0}$ are themselves Markov chains, having jump rates $q(\xi, \xi')$ and $g_t(y, y')$ respectively. Set

$$\tilde{T}_0 = \inf\{t \geqslant 0 : Y_t \neq y_t\} \wedge t_0,$$

then $\tilde{T}_0 > 0$ and, for $t < t_0$, the hazard rate for $\tilde{T}_0$ is given by $\rho(t, X_{t-}, Y_{t-})$, where

$$\rho(t, \xi, y) = \sum_{y' \neq y} |\gamma(\xi, y') - g_t(y, y')|.$$

Thus, there is an exponential random variable $E$ of parameter 1 such that, on $\{\tilde{T}_0 < t_0\}$,

$$E = \int_0^{\tilde{T}_0} \rho(t, X_t, Y_t) dt.$$

On $\Omega_0 \cap \Omega_1 \cap \Omega_2' \cap \Omega_3$ we know that $\sup_{t \leqslant t_0} \|\boldsymbol{X}_t - x_t\| \leqslant \varepsilon$ so, if also $\tilde{T}_0 < \tau$, then $\tilde{T}_0 \leqslant T_0$ and so

$$\int_0^{\tilde{T}_0} \rho(t, X_t, Y_t) dt \leqslant (G + \kappa) t_0.$$

Hence $\mathbb{P}(\tilde{T}_0 < \tau$ and $\Omega_0 \cap \Omega_1 \cap \Omega_2' \cap \Omega_3) \leqslant \mathbb{P}(E \leqslant (G+\kappa)t_0) \leqslant (G+\kappa)t_0$, which combines with our earlier estimates to give the desired result. □

There are at least two places where the basic argument used throughout this section is wasteful and where, with extra effort, better estimates could be obtained. First, we have treated the coordinate functions symmetrically; it may be that a rescaling of some coordinate functions would have the effect of equalizing the noise in each direction. This will tend to improve the estimates. Second, Gronwall's lemma is a blunt instrument. A better idea of how the perturbations introduced by the noise actually propagate is provided by differentiating the solution flow to the differential equation. Sometimes it is possible to show that, rather than growing exponentially, the effect of perturbations actually decays with time. These refinements are particularly relevant, respectively, when the dimension $d$ is large, and when the time horizon $t_0$ is large. We do not pursue them further here.

## 5. Stochastic epidemic

We discuss this well known model, see for example [3], to show in a simple context how the estimates of the preceding section lead quickly to useful asymptotic results. The stochastic epidemic in a population of size $N$ is a Markov chain $X = (X_t)_{t \geqslant 0}$ whose state-space $S$ is the set of pairs $\xi = (\xi^1, \xi^2)$ of non-negative integers with $\xi^1 + \xi^2 \leqslant N$. The non-zero jump rates, for distinct $\xi, \xi' \in S$, are given by

$$q(\xi, \xi') = \begin{cases} \lambda \xi^1 \xi^2 / N, & \text{if } \xi' = \xi + (-1, 1), \\ \mu \xi^2, & \text{if } \xi' = \xi + (0, -1). \end{cases}$$

Here $\lambda$ and $\mu$ are positive parameters, having the interpretation of infection and removal rates, respectively. Write $X_t = (\xi_t^1, \xi_t^2)$. Then $\xi_t^1$ represents the number of susceptible individuals at time $t$ and $\xi_t^2$ the number of infective individuals. Suppose that initially a proportion $p \in (0, 1)$ of the population is infective, the rest being susceptible. Thus $X_0 = (N(1-p), Np)$. The choice of jump rates arises from the modelling assumption that each susceptible individual encounters randomly other members of the population, according to a Poisson process and becomes infective on first meeting an infective individual; then infectives are removed at an exponential rate $\mu$. By a linear change of timescale we can reduce to the case $\mu = 1$, so we shall assume that $\mu = 1$ from this point on.

### *5.1. Convergence to a limit differential equation*

Define $\boldsymbol{x} : S \to \mathbb{R}^2$ by $\boldsymbol{x}(\xi) = \xi/N$ and set $\boldsymbol{X}_t = \boldsymbol{x}(X_t)$. Then the drift vector is given by $\beta(\xi) = b(\boldsymbol{x}(\xi))$, where $b(x) = (-\lambda x^1 x^2, \lambda x^1 x^2 - x^2)$ and $\phi(\xi, \theta) = \sigma_\theta(1/N)(\lambda \xi^1 \xi^2/N + \xi^2)$. Take $U = [0,1]^2$ and set $x_0 = (1-p, p)$. The differential equation $\dot{x}_t = b(x_t)$, which is written in coordinates as

$$\dot{x}_t^1 = -\lambda x_t^1 x_t^2, \quad \dot{x}_t^2 = \lambda x_t^1 x_t^2 - x_t^2,$$

has a unique solution $(x_t)_{t \geq 0}$, starting from $x_0$, which stays in $U$ for all time. Note that $\boldsymbol{x}(S) \subseteq U$, so condition (2) holds for any $\varepsilon > 0$ and $t_0$. The Lipschitz constant for $b$ on $U$ is given by $K = \lambda + \lambda \vee 1$. Set $A = (1+\lambda)e/N$ and take $\delta = e^{-Kt_0}\varepsilon/3$ and $\theta = \delta/(At_0)$, as in Section 4. Let us assume that $\varepsilon \leq t_0$, then $\theta \leq N$, so $\sigma_\theta(1/N) \leq \frac{1}{2}(\theta/N)^2 e$ (as in Footnote 4) and so

$$\int_0^{T \wedge t_0} \phi(X_t, \theta) dt \leq N\sigma_\theta(1/N)(\lambda+1)t_0 \leq \tfrac{1}{2}\theta^2 A t_0.$$

Hence, in this example, $\Omega_0 = \Omega_1 = \Omega_2 = \Omega$ and from Theorem 4.2 we obtain the estimate

$$\mathbb{P}\left(\sup_{t \leq t_0} \|\boldsymbol{X}_t - x_t\| > \varepsilon \right) \leq 4e^{-N\varepsilon^2/C}, \tag{3}$$

where $C = 18(\lambda+1)t_0 e^{2Kt_0+1}$. Figure 2 illustrates a realization of the process alongside the solution of the differential equation.

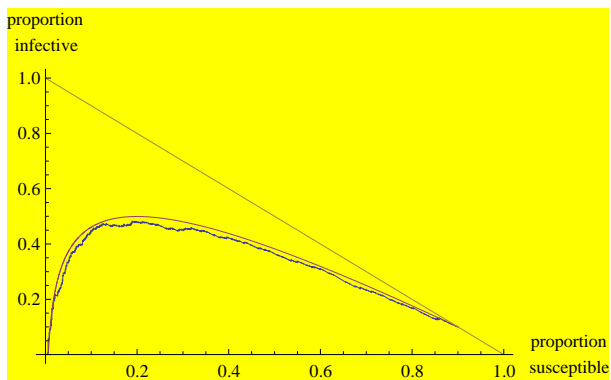

FIG 2. *The graphic shows the proportions of susceptible and infective individuals in a population of 1000, of which initially 900 are susceptible and 100 are infective. The parameter values are $\lambda = 5$ and $\mu = 1$. One realization of the Markov chain, and the solution of the differential equation, are shown at $1 : 1000$ scale.*

### *5.2. Convergence of the terminal value*

The estimate just obtained, whilst giving strong control of error probabilities as $N$ becomes large, behaves rather poorly as a function of $t_0$. This is because we have used the crude device of Gronwall's lemma rather than paying closer attention to the stability properties of the differential equation $\dot{x}_t = b(x_t)$. In particular, the estimate is useless if we want to predict the final size of the epidemic, that is to say, the proportion of the population which is eventually infected, given by $X_\infty^3 = \lim_{t \to \infty} X_t^3$, where $X_t^3 = 1 - X_t^1 - X_t^2$. However, we can obtain an estimate on $X_\infty^3$ by the following modified approach. Let us change the non-zero jump rates by setting $\tilde{q}(\xi, \xi') = q(\xi, \xi')/\xi^2$, for $\xi = (\xi^1, \xi^2)$, to obtain a new process $(\tilde{X}_t)_{t \geq 0}$. Since we have changed only the time-scale, the final values $X_\infty^3$ and $\tilde{X}_\infty^3$ have the same distribution. We can now re-run the analysis, just done for $X$, to $\tilde{X}$. Using obvious notation, we have $\tilde{b}(x) = (-\lambda x^1, \lambda x^1 - 1)$ and $\tilde{\phi}(\xi, \theta) = \sigma_\theta(1/N)(\lambda \xi^1/N + 1)$. We now take $U = (0, 1]^2$. The Lipschitz constant $K$ is unchanged. We make the obvious extension of $\tilde{b}$ to $\mathbb{R}^2$. By explicit solution of the differential equation, we see that $(\tilde{x}_t)_{t \geq 0}$ leaves $U$ at time $\tau$, with $\tilde{x}_\tau^3 = 1 - \tilde{x}_\tau^1 - \tilde{x}_\tau^2 = \tau$, where $\tau$ is the unique root of the equation

$$\tau + (1-p)e^{-\lambda \tau} = 1.$$

Moreover $b^2(x_\tau) = \lambda x_\tau^1 - 1 < 0$, so $b(x\tau)$ is not tangent to the boundary, and so $\varepsilon + \rho(\varepsilon) \leq C\varepsilon$ for all $\varepsilon \in (0, 1]$ for some $C < \infty$ depending only on $\lambda$ and $p$. We can therefore choose $t_0 > \tau$ and apply Theorem 4.3 to obtain, for a constant $C < \infty$ of the same dependence, for all $\varepsilon \in (0, 1]$ and all $N$,

$$\mathbb{P}\left(|X_\infty^3 - \tau| > \varepsilon\right) \leq 4e^{-N\varepsilon^2/C}.$$

### 5.3. Limiting behaviour of individuals

We finally give an alternative analysis which yields a more detailed picture. Consider a Markov chain $\tilde{X} = (\tilde{X}_t)_{t \geq 0}$ with state-space $\tilde{S}$ consisting of $N$-vectors $\eta = (\eta^1, \ldots, \eta^N)$ with $\eta^j \in \{1, 2, 3\}$ for all $j$. Each component of $\eta$ represents the state of an individual member of the population, state 1 corresponding to susceptible, state 2 to infective, and 3 to removed. The non-zero jump rates, for distinct $\eta, \eta' \in \tilde{S}$, are given by

$$q(\eta, \eta') = \begin{cases} \lambda \xi^2(\eta)/N, & \text{if } \eta' = \eta + e_j \text{ for some } j \text{ with } \eta^j = 1, \\ 1, & \text{if } \eta' = \eta + e_j \text{ for some } j \text{ with } \eta^j = 2. \end{cases}$$

Here $\xi^i(\eta) = |\{j : \eta^j = i\}|$, $i = 1, 2$, and $e_j = e_j^{(N)} = (0, \ldots, 1, \ldots, 0)$ is the elementary $N$-vector with a 1 in the $j$th position. Set $X_t = \xi(\tilde{X}_t)$. Then $X = (X_t)_{t \geq 0}$ is the stochastic epidemic considered above. Define $\boldsymbol{x} : \tilde{S} \to \mathbb{R}^2$ by $\boldsymbol{x}(\eta) = \boldsymbol{x}(\xi(\eta))$. Then $\boldsymbol{x}(\tilde{X}_t) = \boldsymbol{x}(X_t) = \boldsymbol{X}_t$, which we already know to remain close to $x_t$ with high probability when $N$ is large.

We can now describe the limiting behaviour of individual members of the population. Fix $k \in \{1, \ldots, N\}$ and set $I = \{1, 2, 3\}^k$. Define $y : \tilde{S} \to I$ by $y(\eta) = (\eta^1, \ldots, \eta^k)$ and set $Y_t = y(\tilde{X}_t)$. We seek to apply Theorem 4.4. Define for $x \in U$ and $n, n' \in \{1, 2, 3\}$

$$g_0(x, n, n') = \begin{cases} \lambda x^2, & \text{if } n = 1 \text{ and } n' = 2, \\ 1, & \text{if } n = 2 \text{ and } n' = 3, \\ 0, & \text{otherwise,} \end{cases}$$

and, for $y, y' \in I$, set $g(x, y, y') = \sum_{j=1}^{k} g_0(x, y^j, y'^j)$. Then the jump rates for $Y$ are given by $\gamma(\eta, y) = g(\boldsymbol{x}(\eta), y(\eta), y)$, so we can take $G = 0$ and $\Omega_3 = \Omega$, and it is straightforward to check that, if $I_0 = I$, then $\kappa = k\lambda\varepsilon$. Hence there is a time-inhomogeneous Markov chain $(y_t)_{t \geqslant 0}$ with state-space $I$ and jump rates $g_t(y, y') = g(x_t, y, y')$, $y, y' \in I$, such that

$$\mathbb{P}(Y_t \neq y_t \text{ for some } t \leqslant t_0) \leqslant k\lambda\varepsilon t_0 + 4e^{-N\varepsilon^2/C}$$

A roughly optimal choice of $\varepsilon$ is $\sqrt{C \log N/N}$, giving a constant $C' < \infty$, depending only on $\lambda$ and $t_0$, such that

$$\mathbb{P}(Y_t \neq y_t \text{ for some } t \leqslant t_0) \leqslant C'k\sqrt{\log N/N}$$

for all sufficiently large $N$. Note that the components of $(y_t)_{t \geqslant 0}$ are independent.

## 6. Population processes

The modelling of population dynamics, involving a number of interacting species, is an important application of Markov chains. A simple example of this was already discussed in Section 5. We propose now to consider another example, of a model which has been used for the growth of a virus in a cell. Our primary aim here is to show how to deal with a Markov chain where some components, the *slow variables*, can be approximated by the solution to a differential equation but others, the *fast variables*, instead oscillate rapidly and randomly. Specifically, by a non-standard choice of coordinate functions, we can obtain an approximation for the slow variables, with computable error probabilities.

A population process is a Markov chain $X = (X_t)_{t \geqslant 0}$, where the state $X_t = (\xi_t^1, \ldots, \xi_t^n)$ describes the number of individuals in each of $n$ species at time $t$; the dynamics are specified by a choice of rates $\lambda_{\varepsilon, \varepsilon'}$ for each of the possible *reactions* $(\varepsilon, \varepsilon')$, where $\varepsilon, \varepsilon' \in (\mathbb{Z}^+)^n$; then, independently over reactions, $X$ makes jumps of size $\varepsilon' - \varepsilon$ at rate

$$\lambda_{\varepsilon, \varepsilon'} \prod_{i=1}^{n} \binom{\xi_i}{\varepsilon_i}.$$

The sort of analysis done below can be adapted to many other population processes.

### 6.1. Analysis of a model for viral replication and growth

We learned of this model from the paper [2], which contains further references on the scientific background. There are three species $G, T$ and $P$ which represent, respectively, the genome, template, and structural protein of a virus. We denote by $\xi^1, \xi^2, \xi^3$ the respective numbers of molecules of each type. There are six reactions, forming a process which may lead from a single virus genome to a sustained population of all three species and to the production of the virus. We write the reactions as follows:

$$G \xrightarrow{\lambda} T, \quad T \xrightarrow{R/\alpha} \emptyset, \quad T \xrightarrow{R} T + G,$$

$$T \xrightarrow{RN} T + P, \quad P \xrightarrow{R/\mu} \emptyset, \quad G + P \xrightarrow{\nu/N} \emptyset.$$

Here, $\alpha > 1, R \geqslant 1, N \geqslant 1$ and $\lambda, \mu, \nu > 0$ are given parameters and, for example, the third reaction corresponds, in the general notation used above, to the case $\varepsilon = (0, 1, 0)$ and $\varepsilon' = (1, 1, 0)$, with $\lambda_{\varepsilon, \varepsilon'} = R$, whereas the final reaction, which causes jumps of size $(-1, 0, -1)$, occurs at a total rate of $\nu \xi^1 \xi^3 / N$. We have omitted some scientific details which are irrelevant to the mathematics, and have written $\emptyset$ when the reaction produces none of the three species in the model. In fact it is the final reaction $G + P$ which gives rise to the virus itself. In the case of scientific interest, $\alpha, \lambda, \mu, \nu$ are of order 1, but $R, N$ are large. We therefore seek an approximation which is good in this regime.

As a first step to understanding this process, we note that, for as long as the number of templates $\xi_t^2$ remains of order 1, the rate of production of genomes is of order $R$. On the other hand, for as long as the number of genomes $\xi_t^1$ is bounded by $xR$, for some $x > 0$, the number of templates can be dominated[7] by a $M/M/\infty$ queue, $(Y_t)_{t \geqslant 0}$ with arrival rate $\lambda x R$ and service rate $R/\alpha$. The stationary distribution for $(Y_t)_{t \geqslant 0}$ is Poisson of parameter $\lambda x \alpha$, which suggests that, for reasonable initial conditions at least, $\xi_t^2$ does remain of order 1, but oscillates rapidly, on a time-scale of order $1/R$. The number of proteins $\xi_t^3$ evolves as an $M/M/\infty$ queue, with time-dependent arrival rate $RN\xi_t^2$ and service rate $R/\mu + \nu \xi_t^1 / N$. This suggests that $\xi_t^3 / N$ will track closely a function of the rapidly oscillating process $\xi_t^2$.

The only hope for a differential equation approximation would thus appear to be the genome process $(\xi_t^1)_{t \geqslant 0}$. The obvious choice of coordinate map $\boldsymbol{x}(\xi) = \xi^1 / R$ gives as drift

$$\beta(\xi) = \sum_{\xi' \neq \xi} (\boldsymbol{x}(\xi') - \boldsymbol{x}(\xi)) q(\xi, \xi') = -\lambda \frac{\xi^1}{R} + \xi^2 - \nu \frac{\xi^1}{R} \frac{\xi^3}{N},$$

which we cannot approximate by a function of $\boldsymbol{x}(\xi)$ unless the second and third terms become negligible. In fact they do not, so this choice fails. The problem is that the drift of $\xi_t^1$ is significantly dependent on the fast variables $\xi_t^2$ and $\xi_t^3$. To overcome this, we can attempt to compensate the coordinate process so that it takes account of this effect. We seek to find a function $\boldsymbol{x}$ on the state-space $S = (\mathbb{Z}^+)^3$ of the form

$$\boldsymbol{x}(\xi) = \frac{\xi^1}{R} + \chi(\xi),$$

---

[7]The obvious additivity property for arrival rates of $M/M/\infty$ queues having a common service rate extends to the case of previsible arrival rates. A good way to see this is by constructing all queues from a single Poisson random measure

where $\chi$ is a small correction, chosen so that the drift vector $\beta(\xi)$ has the form

$$\beta(\xi) = b(\boldsymbol{x}(\xi)) + \frac{\Delta(\xi)}{R},$$

where again $\Delta(\xi)/R$ is small when $R$ is large. Small here refers to a typical evaluation on the process, where we recall that we expect $\xi_t^1/R$, $\xi_t^2$ and $\xi_t^3/N$ to be of order 1. It is reasonable to search for a function $\chi$ which is affine in $\xi^1$ and linear in $(\xi^2, \xi^3)$. After some straightforward calculations, we find that

$$\chi(\xi) = \frac{1}{R}\left(\alpha\xi^2 - \mu\nu\frac{\xi^1}{R}\frac{\xi^3}{N} - \alpha\mu\nu\frac{\xi^1}{R}\xi^2\right)$$

has the desired property, with

$$b(x) = \lambda(\alpha - 1)x - \lambda\alpha\mu\nu x^2$$

and

$$\Delta = \lambda\mu\nu\frac{\xi^1}{R}\frac{\xi^3}{N} + \alpha\lambda\mu\nu(\xi^2 + 1)\frac{\xi^1}{R} - \mu\nu\xi^2\frac{\xi^3}{N} - \alpha\mu\nu(\xi^2)^2 + \alpha\mu\nu^2\frac{\xi^1}{R}\xi^2\frac{\xi^3}{N}$$
$$+ \mu\nu^2\frac{\xi^1}{R}\frac{\xi^3}{N}\frac{(\xi^1 + \xi^3 - 1)}{N} - \lambda(\alpha - 1)R\chi(\xi) + \lambda\alpha\mu\nu\left(2R\chi(\xi)\frac{\xi^1}{R} + R\chi(\xi)^2\right).$$

The limit differential equation

$$\dot{x}_t = \lambda(\alpha - 1)x_t - \lambda\alpha\mu\nu x_t^2$$

has a unique positive fixed point $x_\infty = (\alpha - 1)/(\alpha\mu\nu)$. Fix $x_0 \in [0, x_\infty]$ and take as initial state $X_0 = (Rx_0, 0, 0)$.

**Theorem 6.1.** *For all $t_0 \in [1, \infty)$, there is a constant $C < \infty$, depending only on $\alpha, \lambda, \mu, \nu, t_0$ with the following property. For all $\varepsilon \in (0, 1]$ there is a constant $R_0 < \infty$, depending only on $\alpha, \lambda, \mu, \nu, t_0$ and $\varepsilon$ such that, for all $R \geqslant R_0$ and $N \geqslant R$, we have*

$$\mathbb{P}\left(\sup_{t \leqslant t_0}\left|\frac{\xi_t^1}{R} - x_t\right| > \varepsilon\right) \leqslant e^{-R\varepsilon^2/C}.$$

*Proof.* We shall write $C$ for a finite constant depending only on $\alpha, \lambda, \mu, \nu, t_0$, whose value may vary from line to line, adding a subscript when we wish to refer to a particular value at a later point. Fix constants $a \geqslant 1$, $\gamma > 0, \Gamma \geqslant 1$, with $(\alpha + 1)(\mu\nu + 1)\gamma \leqslant 1/2$, to be determined later, and set $A = a/R$. Take $U = [0, x_\infty + 1]$. As in Section 4, let us write $K$ for the Lipschitz constant of $b$ on $U$, and set $\boldsymbol{X}_t = \boldsymbol{x}(X_t)$, $\delta = \varepsilon e^{-Kt_0}/3$, $\theta = \delta/(At_0)$ and $T = \inf\{t \geqslant 0 : \boldsymbol{X}_t \notin U\}$. Since $0 \leqslant x_t \leqslant x_\infty$ for all $t$ and since $\varepsilon \leqslant 1$, condition (2) holds.

Consider the events

$$\Omega_4 = \left\{\sup_{0 \leqslant t \leqslant T \wedge t_0} \xi_t^2 \leqslant \gamma R \quad \text{and} \quad \sup_{0 \leqslant t \leqslant T \wedge t_0} \xi_t^3 \leqslant \gamma RN\right\}$$

and
$$\Omega_5 = \left\{ \int_0^{T \wedge t_0} \xi_t^2 dt \leqslant \Gamma \quad \text{and} \quad \int_0^{T \wedge t_0} \xi_t^3 dt \leqslant \Gamma N \right\}.$$

We refer to Subsection 4.2 for the definition of the events $\Omega_0, \Omega_1, \Omega_2$. We now show that, for suitable choices of the constants $a, \gamma, \Gamma$, we have $\Omega_4 \cap \Omega_5 \subseteq \Omega_0 \cap \Omega_1 \cap \Omega_2$.

For $t \leqslant T \wedge t_0$, on $\Omega_4$, we have

$$\boldsymbol{X}_t = \frac{\xi_t^1}{R}\left(1 - \alpha\mu\nu\frac{\xi_t^2}{R} - \mu\nu\frac{\xi_t^3}{RN}\right) + \frac{\alpha\xi_t^2}{R} \geqslant \frac{\xi_t^1}{2R}, \tag{4}$$

so $\xi_t^1/R \leqslant 2(x_\infty + 1)$.

On $\Omega_4$ (without using the assumptions $\xi_0^2 = \xi_0^3 = 0$), we have

$$|\boldsymbol{X}_0 - x_0| = |\chi(X_0)| = \frac{1}{R}\left|\alpha\xi_0^2 - \mu\nu\frac{\xi_0^1}{R}\frac{\xi_0^3}{N} - \alpha\mu\nu\frac{\xi_0^1}{R}\xi_0^2\right| \leqslant C_0 \gamma,$$

so, provided that $C_0 \gamma \leqslant \delta$, we have $\Omega_4 \subseteq \Omega_0$.

On $\Omega_4 \cap \Omega_5$, we have

$$\int_0^{T \wedge t_0} |\beta(X_t) - b(\boldsymbol{x}(X_t))| dt = \frac{1}{R}\int_0^{T \wedge t_0} |\Delta(X_t)| dt$$
$$\leqslant \frac{C}{R}\int_0^{T \wedge t_0}\left(1 + |\xi_t^2|^2 + \left|\frac{\xi_t^3}{N}\right|^2\right)dt \leqslant \frac{C_1}{R}(1 + R\gamma\Gamma).$$

So, provided that $C_1(\gamma\Gamma + 1/R) \leqslant \delta$, we have $\Omega_4 \cap \Omega_5 \subseteq \Omega_1$.

For $\xi \in S$ with $\xi^1 \leqslant (x_\infty + 1)R$, $\xi^2 \leqslant \gamma R$ and $\xi^3 \leqslant \gamma RN$, and for any $\xi' \in S$ with $q(\xi, \xi') > 0$, we have $|\xi'^i - \xi^i| \leqslant 1$ and hence

$$|\boldsymbol{x}(\xi') - \boldsymbol{x}(\xi)| = \left|\frac{\xi'^1 - \xi^1}{R} + \chi(\xi') - \chi(\xi)\right| \leqslant \frac{C}{R},$$

and indeed, for $\xi' = \xi \pm (0, 0, 1)$ we have $|\boldsymbol{x}(\xi') - \boldsymbol{x}(\xi)| \leqslant C/(RN)$, so, using

$$\sigma_\theta(\boldsymbol{x}(\xi') - \boldsymbol{x}(\xi)) \leqslant \frac{1}{2}\frac{C\theta^2}{R}e^{C\theta/R}|\boldsymbol{x}(\xi') - \boldsymbol{x}(\xi)|,$$

we obtain, after some straightforward estimation,

$$\phi(\xi, \theta) \leqslant \frac{1}{2}\frac{C\theta^2}{R}e^{C\theta/R}\sum_{\xi' \neq \xi}|\boldsymbol{x}(\xi') - \boldsymbol{x}(\xi)|q(\xi, \xi') \leqslant \frac{1}{2}\frac{C\theta^2}{R}e^{C\theta/R}\left(1 + \xi^2 + \frac{\xi^3}{N}\right).$$

So, on $\Omega_4 \cap \Omega_5$, we have

$$\int_0^{T \wedge t_0} \phi(X_t, \theta) dt \leqslant \frac{1}{2}\frac{C\theta^2}{R}e^{C\theta/R}\int_0^{T \wedge t_0}\left(1 + \xi_t^2 + \frac{\xi_t^3}{N}\right)dt \leqslant \frac{1}{2}\frac{C_2\Gamma\theta^2}{R}e^{C_2\theta/R}.$$

So, provided $a \geqslant C_2\Gamma e$, we have $C_2\theta/R \leqslant 1$ and so $\Omega_4 \cap \Omega_5 \subseteq \Omega_2$.

From equation (4) we obtain $|\xi_t^1/R - \boldsymbol{X}_t| \leqslant C_3\gamma$. Let us choose then $a, \gamma, \Gamma$ and $R$ so that $C_0\gamma \leqslant \delta$, $C_1(\gamma\Gamma + 1/R) \leqslant \delta$, $a \geqslant C_2\Gamma e$ and $C_3\gamma \leqslant \varepsilon$. Then

$$\mathbb{P}\left(\sup_{t \leqslant t_0} \left|\frac{\xi_t^1}{R} - x_t\right| > 2\varepsilon\right) \leqslant \mathbb{P}\left(\sup_{t \leqslant t_0} |\boldsymbol{X}_t - x_t| > \varepsilon\right) + \mathbb{P}(\Omega_4^c).$$

On the other hand, $\Omega_4 \cap \Omega_5 \subseteq \Omega_0 \cap \Omega_1 \cap \Omega_2$ and, by Theorem 4.2, we have

$$\mathbb{P}\left(\sup_{t \leqslant t_0} |\boldsymbol{X}_t - x_t| > \varepsilon\right) \leqslant 2e^{-\delta^2/(2At_0)} + \mathbb{P}(\Omega_0^c \cup \Omega_1^c \cup \Omega_2^c).$$

Since $2e^{-\delta^2/(2At_0)} = 2e^{-R\varepsilon^2/C}$, where $C = 18t_0e^{2Kt_0}$, we can now complete the proof by showing that, for suitable $a, \gamma, \Gamma$ and $R_0$, for all $R \geqslant R_0$, we have $\mathbb{P}(\Omega_4^c) \leqslant e^{-R}$ and $\mathbb{P}(\Omega_5^c) \leqslant e^{-R}$.

We can dominate the processes $(\xi_t^2)_{t \geqslant 0}$ and $(\xi_t^2)_{t \geqslant 0}$, up to $T$, by a pair of processes $Y = (Y_t)_{t \geqslant 0}$ and $Z = (Z_t)_{t \geqslant 0}$, respectively, where $Y$ is an $M/M/\infty$ queue with arrival rate $2\lambda(x_\infty + 1)R$ and service rate $R/\alpha$, starting from $\xi_0^2 = 0$, and where, conditional on $Y$, $Z$ is an $M/M/\infty$ queue with arrival rate $RNY_t$ and service rate $R/\mu$, starting from $\xi_0^3 = 0$. We now use the estimates (5) and (6), to be derived in the next subsection. For $\Gamma$ sufficiently large, using the estimate (6), we have $\mathbb{P}(\Omega_5^c) \leqslant e^{-R}$ for all sufficiently large $R$. Fix such a $\Gamma$ and choose $a, R$ sufficiently large and $\gamma$ sufficiently small to satisfy the above constraints. Finally, using the estimate (5), $\mathbb{P}(\Omega_4^c) \leqslant e^{-R}$, for all sufficiently large $R$. $\square$

The initial state $(Rx_0, 0, 0)$ was chosen to simplify the presentation and is not realistic. However, an examination of the proof shows that, for some constant $\gamma > 0$, depending only on $\alpha, \lambda, \mu, \nu, \varepsilon$, the same conclusion can be drawn for any initial state $(Rx_0, \xi_0^2, \xi_0^3)$ with $x_0 \leqslant x_\infty$, $\xi_0^2 \leqslant R\gamma$ and $\xi_0^3 \leqslant RN\gamma$. Since typical values of the fast variables $\xi_t^2$ and $\xi_t^3$ are of order 1 and $N$ respectively, this is more realistic. Although we are free to take an initial state $(1, 0, 0)$, the action of interest in this case occurs at a time of order $\log R$, so is not covered by our result. Instead, there is a branching process approximation for the number of genomes, valid until it reaches $Rx_0$, for small $x_0$. Our estimate can be applied to the evolution of the process from that time on. See [2] for more details of the branching process approximation.

### *6.2. Some estimates for the $M/M/\infty$ queue*

We now derive the fast variable bounds used in the proof of Theorem 6.1. They are based on the following two estimates for the $M/M/\infty$ queue.

**Proposition 6.2.** *Let $(X_t)_{t \geqslant 0}$ be an $M/M/\infty$ queue starting from $x_0$, with arrival rate $\lambda$ and service rate $\mu$. Then, for all $t \geqslant 1/\mu$ and all $a \geqslant 3\lambda e^2/\mu$,*

$$\mathbb{P}\left(\sup_{0 \leqslant s \leqslant t} X_s \geqslant x_0 + \log(\mu t) + a\right) \leqslant \exp\left\{-a\log\left(\frac{\mu a}{3\lambda e}\right)\right\}.$$

*Proof.* By rescaling time, we reduce to the case where $\mu = 1$. Also, we have $X_t \leqslant x_0 + Y_t$, where $(Y_t)_{t \geqslant 0}$ is an $M/M/\infty$ queue starting from 0, with the same arrival rate and service rate. Thus we are reduced to the case where $x_0 = 0$.

Choose $n \in \mathbb{N}$ so that $\delta = t/n \in [1,2)$ and note that, for $k = 0, 1, \ldots, n-1$, we have
$$\sup_{k\delta \leqslant s \leqslant (k+1)\delta} X_s \leqslant X_{k\delta} + A_{k+1} \leqslant Y_k + A_{k+1},$$
where $A_{k+1}$ is the number of arrivals in $(k\delta, (k+1)\delta]$ and where $Y_k$ is a Poisson random variable of parameter $\lambda$, independent of $A_{k+1}$. By the usual Poisson tail estimate[8], for all $x \geqslant 0$,
$$\mathbb{P}\left(\sup_{k\delta \leqslant s \leqslant (k+1)\delta} X_s \geqslant x\right) \leqslant \exp\left\{-x \log\left(\frac{x}{\lambda(1+\delta)e}\right)\right\}.$$

Hence, for $t \geqslant 1$ and $a \geqslant 3\lambda e^2$,
$$\mathbb{P}\left(\sup_{0 \leqslant s \leqslant t} X_s \geqslant \log t + a\right)$$
$$\leqslant n \exp\left\{-(a + \log t)\log\left(\frac{a}{3\lambda e}\right)\right\} \leqslant \exp\left\{-a \log\left(\frac{a}{3\lambda e}\right)\right\}.$$
□

**Proposition 6.3.** *Let $(X_t)_{t \geqslant 0}$ be an $M/M/\infty$ queue starting from $x_0$, with time-dependent arrival rate $\lambda_t$ and service rate $\mu$. Then, for all $t \geqslant 0$ and all $\theta \in [0, \mu)$,*
$$\mathbb{E}\left(\exp\left\{\theta \int_0^t X_s ds\right\}\right) \leqslant \left(\frac{\mu}{\mu - \theta}\right)^{x_0} \exp\left\{\frac{\theta}{\mu - \theta} \int_0^t \lambda_s ds\right\}$$

*Proof.* By rescaling time, we reduce to the case where $\mu = 1$. Consider first the case where $\lambda_t \equiv 0$. Then
$$\int_0^\infty X_s ds = S_1 + \cdots + S_{x_0}$$
where $S_n$ is the service time of the $n$th customer present at time 0. The result follows by an elementary calculation.

Consider next the case where $x_0 = 0$. We can express $X_t$ in terms of a Poisson random measure $m$ on $[0, \infty) \times [0,1]$ of intensity $\lambda_t dt du$, thus
$$X_t = \int_0^t \int_0^1 1_{\{u \leqslant e^{-(t-s)}\}} m(ds, du).$$
Then, by Fubini,
$$\int_0^t X_s ds \leqslant \int_0^t \int_0^1 \log\left(\frac{1}{u}\right) m(ds, du).$$

---
[8] For a Poisson random variable $X$ of parameter $\lambda$, we have $\mathbb{P}(X \geqslant x) \leqslant \exp\{-x \log(\frac{x}{\lambda e})\}$, for all $x \geqslant 0$.

By Campbell's formula,

$$\mathbb{E}\left(\exp\left\{\theta \int_0^t \int_0^1 \log\left(\frac{1}{u}\right) m(ds,du)\right\}\right) = \exp\left\{\frac{\theta}{1-\theta}\int_0^t \lambda_s ds\right\}.$$

The result for $x_0 = 0$ follows. The general case now follows by independence. □

Now let $Y = (Y_t)_{t \geqslant 0}$ be an $M/M/\infty$ queue with arrival rate $\lambda R$ and service rate $R/\alpha$, starting from $Ry$, and, conditional on $Y$, let $(Z_t)_{t \geqslant 0}$ be an $M/M/\infty$ queue with arrival rate $RNY_t$ and service rate $R/\mu$, starting from $RNz$. By Proposition 6.2, for any $t_0 \geqslant 0$ and $\gamma > 0$, we can find $R_0 < \infty$, depending only on $\alpha, \gamma, \lambda, \mu$ and $t_0$, such that, for all $R_0 \geqslant R$ and $N \geqslant 1$,

$$\mathbb{P}\left(\sup_{0 \leqslant t \leqslant t_0} Y_t \geqslant (\gamma + y)R \text{ or } \sup_{0 \leqslant t \leqslant t_0} Z_t \geqslant (\gamma + y + z)RN\right) \leqslant e^{-R}. \quad (5)$$

On the other hand, using Proposition 6.3, for $\theta\alpha \leqslant 1/2$,

$$\mathbb{E}\left(\exp\left\{R\theta \int_0^{t_0} Y_t dt\right\}\right) \leqslant \left(\frac{1}{1-\theta\alpha}\right)^y \exp\left\{\frac{\theta\alpha\lambda R t_0}{1-\theta\alpha}\right\} \leqslant e^{CR(1+y)},$$

for a constant $C$ depending on $\alpha, \lambda, t_0$. Then, by conditioning first on $Y$, we obtain, for all $\theta(\mu+1)\alpha \leqslant 1/2$,

$$\mathbb{E}\left(\exp\left\{\frac{R\theta}{N}\int_0^{t_0} Z_t dt\right\}\right)$$
$$\leqslant \left(\frac{1}{1-\theta\mu/N}\right)^{RNz}\mathbb{E}\left(\exp\left\{\frac{R\mu\theta}{1-\theta\mu/N}\int_0^{t_0} Y_t dt\right\}\right) \leqslant e^{C(1+y+z)},$$

where $C$ depends on $\alpha, \lambda, \mu, t_0$. Thus, we obtain constants $\Gamma, R_0 < \infty$, depending only on $\alpha, \lambda, \mu, t_0$, such that, for all $R \geqslant R_0$ and $N \geqslant 1$,

$$\mathbb{P}\left(\int_0^{t_0} Y_t dt \geqslant \Gamma(1+y) \text{ or } \int_0^{t_0} Z_t dt \geqslant \Gamma(1+y+z)N\right) \leqslant e^{-R}. \quad (6)$$

## 7. Hypergraph cores

The approximation of Markov chains by differential equations is a powerful tool in probabilistic combinatorics, and in particular in the asymptotic analysis of structures within large random graphs and hypergraphs. It is sometimes possible to find an algorithm, whose progress can be described in terms of a Markov chain, and whose terminal value gives information about the structure of interest. If this Markov chain can be approximated by a differential equation, then this may provide an effective means of computation. We shall describe in detail an implementation of this approach which yields a quantitative description of the $k$-core for a general class of random hypergraphs. Here $k \geqslant 2$ is an integer, which will remain fixed throughout.

### 7.1. Specification of the problem

Let $V$ and $E$ be finite sets. A *hypergraph* with *vertex set* $V$ and *edge-label set* $E$ is a subset $\gamma$ of $V \times E$. Given a hypergraph $\gamma$, define, for $v \in V$ and $e \in E$, sets $\gamma(v) = \gamma \cap (\{v\} \times E)$ and $\gamma(e) = \gamma \cap (V \times \{e\})$. The sets $\gamma(e)$ are the (*hyper*)*edges* of the hypergraph $\gamma$. Figure 3 gives two pictorial representations of a small hypergraph. The *degree* and *weight* functions $\boldsymbol{d}_\gamma : V \to \mathbb{Z}^+$ and $\boldsymbol{w}_\gamma : E \to \mathbb{Z}^+$ of $\gamma$ are given by $\boldsymbol{d}_\gamma(v) = |\gamma(v)|$ and $\boldsymbol{w}_\gamma(e) = |\gamma(e)|$. The *k-core* $\bar{\gamma}$ of $\gamma$ is the largest subset $\bar{\gamma}$ of $\gamma$ such that, for all $v \in V$ and $e \in E$,

$$\boldsymbol{d}_{\bar{\gamma}}(v) \in \{0\} \cup \{k, k+1, \dots\}, \quad \boldsymbol{w}_{\bar{\gamma}}(e) \in \{0, \boldsymbol{w}_\gamma(e)\}.$$

Thus, if we call a *sub-hypergraph* of $\gamma$ any hypergraph obtained by deleting edges from $\gamma$, then $\bar{\gamma}$ is the largest sub-hypergraph of $\gamma$ in which every vertex of non-zero degree has degree at least $k$. It is not hard to see that any algorithm which deletes recursively edges containing at least one vertex of degree less than $k$ terminates at the $k$-core $\bar{\gamma}$. The $k$-core is of interest because it is a measure of the strength of connectivity present in $\gamma$; see [17], [19], [20].

A *frequency vector* is a vector $\boldsymbol{n} = (n_d : d \in \mathbb{Z}^+)$ with $n_d \in \mathbb{Z}^+$ for all $d$. We write $m(\boldsymbol{n}) = \sum_d d n_d$. Given a function $\boldsymbol{d} : V \to \mathbb{Z}^+$, define its frequency vector $\boldsymbol{n}(\boldsymbol{d}) = (n_d(\boldsymbol{d}) : d \in \mathbb{Z}^+)$ by $n_d(\boldsymbol{d}) = |\{v \in V : \boldsymbol{d}(v) = d\}|$, and set $m(\boldsymbol{d}) = m(\boldsymbol{n}(\boldsymbol{d})) = \sum_v \boldsymbol{d}(v)$. The *frequency vectors* of a hypergraph $\gamma$ are then the pair $\boldsymbol{p}(\gamma), \boldsymbol{q}(\gamma)$, where $\boldsymbol{p}(\gamma) = \boldsymbol{n}(\boldsymbol{d}_\gamma)$ and $\boldsymbol{q}(\gamma) = \boldsymbol{n}(\boldsymbol{w}_\gamma)$. Note that $m(\boldsymbol{p}(\gamma))$ is simply the cardinality of $\gamma$, as of course is $m(\boldsymbol{q}(\gamma))$.

The datum for our model is a pair of non-zero frequency vectors $\boldsymbol{p}, \boldsymbol{q}$ with $m(\boldsymbol{p}) = m(\boldsymbol{q}) = m < \infty$. Note that there exists an integer $L \geqslant 2$ such that $p_d = q_w = 0$ for all $d, w \geqslant L + 1$. We assume also that $p_0 = q_0 = 0$. This will result in no essential loss of generality. Fix an integer $N \geqslant 1$. We shall be interested in the limit as $N \to \infty$. Choose sets $V$ and $E$ and functions $\boldsymbol{d} : V \to \mathbb{Z}^+$ and $\boldsymbol{w} : E \to \mathbb{Z}^+$ such that $\boldsymbol{n}(\boldsymbol{d}) = N\boldsymbol{p}$ and $\boldsymbol{n}(\boldsymbol{w}) = N\boldsymbol{q}$. In particular, this implies that $|V| = N \sum_d p_d$ and $|E| = N \sum_w q_w$. Denote by $G(\boldsymbol{d}, \boldsymbol{w})$ the set of hypergraphs on $V \times E$ with degree function $\boldsymbol{d}$ and weight function $\boldsymbol{w}$. Thus

$$G(\boldsymbol{d}, \boldsymbol{w}) = \{\gamma \subseteq V \times E : \boldsymbol{d}_\gamma = \boldsymbol{d}, \boldsymbol{w}_\gamma = \boldsymbol{w}\}$$

and, in particular, all elements of $G(\boldsymbol{d}, \boldsymbol{w})$ have cardinality $Nm$. This set is known to be non-empty for $N$ sufficiently large. Its elements can also be thought of as bipartite graphs on $V \cup E$ with given degrees. We shall be interested in the distribution of the $k$-core $\bar{\Gamma}$ when $\Gamma$ is a hypergraph chosen uniformly at random from $G(\boldsymbol{d}, \boldsymbol{w})$. We write $\Gamma \sim U(\boldsymbol{d}, \boldsymbol{w})$ for short. Set $\bar{\boldsymbol{D}} = \boldsymbol{d}_{\bar{\Gamma}}$ and $\bar{\boldsymbol{W}} = \boldsymbol{w}_{\bar{\Gamma}}$. These are the degree and weight functions of the $k$-core. Define for $d, d', w \geqslant 0$

$$\bar{P}_{d,d'} = |\{v \in V : \boldsymbol{d}(v) = d', \bar{\boldsymbol{D}}(v) = d\}|/N, \quad \bar{Q}_w = |\{e \in E : \bar{\boldsymbol{W}}(e) = w\}|/N. \tag{7}$$

Note that, given $(\bar{P}_{d,d'} : k \leqslant d \leqslant d')$ and $(\bar{Q}_w : w \geqslant 1)$, we can recover the other non-zero frequencies from the equations

$$\bar{P}_{0,d'} = p_{d'} - \sum_{k \leqslant d \leqslant d'} \bar{P}_{d,d'}, \quad \bar{Q}_0 = \sum_{w \geqslant 0} q_w - \sum_{w \geqslant 1} \bar{Q}_w,$$

and, given all these frequencies, the joint distribution of $\Gamma$ and its $k$-core $\bar{\Gamma}$ is otherwise dictated by symmetry[9]. For $\bar{D}$ and $\bar{W}$ are independent and uniformly distributed, subject to the equations (7) and to the constraint $\bar{W}(e) \in \{0, w(e)\}$ for all $e \in E$. Moreover, we shall see that, given $\bar{D}$ and $\bar{W}$, $\bar{\Gamma} \sim U(\bar{D}, \bar{W})$. The problem of characterizing the distribution of the the $k$-core thus reduces to that of understanding the frequencies $(\bar{P}_{d,d'} : k \leqslant d \leqslant d')$ and $(\bar{Q}_w : w \geqslant 1)$.

### 7.2. Branching process approximation

In this subsection we describe an approximation to the local structure of a hypergraph $\Gamma \sim U(d, w)$ on which the later analysis relies, and which is illustrated in Figure 3. We work in a more general set-up than the sequence parametrized by $N$ just described. Fix $L < \infty$ and degree and weight functions $d, w$, with $m(d) = m(w) = m$. We consider the limit $m \to \infty$ subject to $d, w \leqslant L$. Note that this limit applies to the set-up of the preceding subsection, where $m(d) = Nm$ with $m$ fixed and $N \to \infty$. Choose a random vertex $v$ according to the distribution $d/m$ and set $D = d(v)$. Enumerate randomly the subset $\Gamma(v) = \{(v, e_1), \ldots, (v, e_D)\}$ and set $S_i = w(e_i) - 1$, $i = 1, \ldots, D$. For $i = 1, \ldots, D$, enumerate randomly the set of vertices in $\Gamma(e_i)$ which are distinct from $v$, thus

$$\Gamma(e_i) = \{v, v_{i,1}, \ldots, v_{i,S_i}\} \times \{e_i\},$$

and set $L_{i,j} = d(v_{i,j}) - 1$. Write $A$ for the event that the vertices $v_{i,j}$ are all distinct. Thus

$$A = \{(v', e_i), (v', e_j) \in \Gamma \text{ implies } v' = v \text{ or } i = j\}.$$

Let $T$ be a discrete alternating random tree, having types $V, E$, with degree distributions $\tilde{p}, \tilde{q}$ respectively, and having base point $\tilde{v}$ of type $V$. Here $\tilde{p}, \tilde{q}$ are the size-biased distributions obtained from $p = n(d)$ and $q = n(w)$ by $\tilde{p}_d = dp_d/m, \tilde{q}_w = wq_w/m, d, w \geqslant 0$. This may be considered as a branching process starting from the single individual $\tilde{v}$, which has $\tilde{D}$ offspring $\tilde{e}_1, \ldots, \tilde{e}_{\tilde{D}}$ of type $E$, where $\tilde{D}$ has distribution $\tilde{p}$; then all individuals of type $E$ have offspring of type $V$, the numbers of these being independent and having distribution $\sigma$; all individuals of type $V$ have offspring, of type $E$, the numbers of these being

---

[9]For the marginal distribution of the $k$-core, only the frequencies

$$\bar{P}_d = |\{v \in V : \bar{D}(v) = d\}|/N = \sum_{d' \geqslant d} \bar{P}_{d,d'}$$

are relevant, but the asymptotics of $P_d$ turn out to split naturally over $d'$, see (11) below.

independent and, with the exception of $\tilde{v}$, having distribution $\lambda$. Here $\lambda$ and $\sigma$ are given by

$$\lambda_d = (d+1)p_{d+1}/m, \quad \sigma_w = (w+1)q_{w+1}/m, \quad d, w \geqslant 0. \tag{8}$$

For $i = 1, \ldots, \tilde{D}$, write $\tilde{S}_i$ for the number of offspring of $\tilde{e}_i$ and, for $j = 1, \ldots, \tilde{S}_i$, write $\tilde{L}_{i,j}$ for the number of offspring of the $j$th offspring of $\tilde{e}_i$. Then, conditional on $\tilde{D} = d$, the random variables $\tilde{S}_1, \ldots, \tilde{S}_d$ are independent, of distribution $\sigma$, and, further conditioning on $\tilde{S}_i = s_i$ for $i = 1, \ldots, d$, the random variables $\tilde{L}_{i,j}, i = 1, \ldots, d, j = 1, \ldots, s_i$, are independent, of distribution $\lambda$.

It is known (see [15] or, for a more explicit statement, [5]) that *there is a function $\psi_0 : \mathbb{N} \to [0, 1]$, depending only on $L$, with $\psi_0(m) \to 0$ as $m \to \infty$, such that, for all degree and weight functions $\boldsymbol{d}, \boldsymbol{w} \leqslant L$ with $m(\boldsymbol{d}) = m(\boldsymbol{w}) = m$, we have $\mathbb{P}(A) \geqslant 1 - \psi_0(m)$ and there is a coupling of $\Gamma$ and $T$ such that $D = \tilde{D}$ and, with probability exceeding $1 - \psi_0(m)$, we have $S_i = \tilde{S}_i$ for all $i$ and $L_{i,j} = \tilde{L}_{i,j}$ for all $i, j$.*

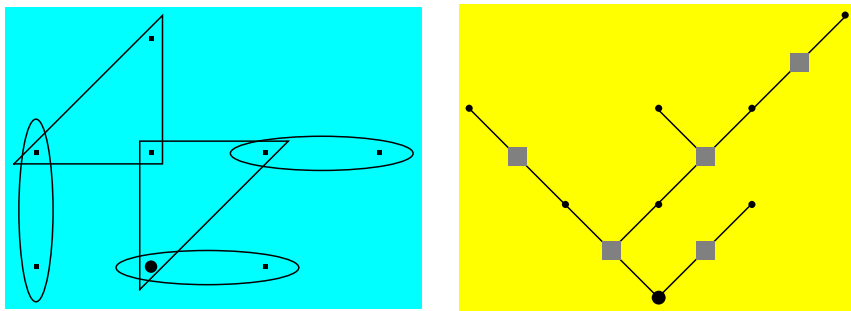

FIG 3. *The left picture shows a hypergraph with eight vertices, three 2-edges, and three 3-edges. An incidence is selected at random, shown by the enlarged vertex, and chosen as root of a branching process, shown as the bottom vertex on the right. The root has two hyperedge offspring, shown as grey squares. One of these has two vertex offspring, and so on.*

The following paragraph presents a heuristic argument which leads quickly to a prediction for the asymptotic frequencies of core degrees and weights, which we shall later verify rigorously, subject to an additional condition. The convergence of $\Gamma$ to $T$, near a randomly chosen vertex, which we expressed in terms of the function $\psi_0$ for the first two steps, in fact holds in a similar sense for any given numbers of steps. The algorithm of deleting, recursively, all edges in $\Gamma$ containing any vertex of degree less than $k$ terminates at the $k$-core $\bar{\Gamma}$. Consider the following analogous algorithm on the branching process: we remove in the first step all individuals of type $E$ having some offspring with fewer than $k - 1$ offspring of its own; then repeat this step infinitely often. Set $g_0 = 1$. For $n \geqslant 0$, write $s_n$ for the probability that, after $n$ steps, a given individual of type $E$ remains in the population, and write $g_{n+1}$ for the probability that, after $n$ steps, a given individual of type $V$ (distinct from $\tilde{v}$) has at least $k - 1$ offspring

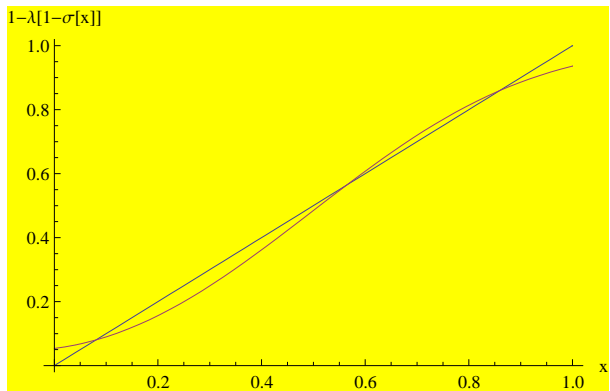

FIG 4. *The function $\phi(x) = 1 - \lambda(1 - \sigma(x))$ is shown for $x \in (0,1)$, where $\lambda(y)$ is a truncated Poisson probability generating function with mean 2.75, and $\sigma(x) = 0.02 + 0.08x + 0.6x^2 + 0.2x^3 + 0.1x^4$. The largest intersection with the line $y = x$ gives the value $g^*$ needed for the 2-core fluid limit.*

remaining. Then, by a standard type of branching process argument,

$$s_n = \sigma(g_n), \quad g_{n+1} = \sum_{j \geqslant k-1} \sum_d \binom{d}{j} \lambda_d \, s_n^j (1-s_n)^{d-j}, \quad n \geqslant 0.$$

We write here $\sigma$, and below $\lambda$, for the probability generating functions

$$\sigma(z) = \sum_w \sigma_w z^w, \quad \lambda(z) = \sum_d \lambda_d z^d. \tag{9}$$

So $g_{n+1} = \phi(g_n)$, where

$$\phi(g) = \sum_{j \geqslant k-1} \sum_d \binom{d}{j} \lambda_d \, \sigma(g)^j (1-\sigma(g))^{d-j}. \tag{10}$$

Note that, in the case $k = 2$, we have the simple formula

$$\phi(g) = 1 - \lambda(1 - \sigma(g)).$$

Since $\phi$ maps $[0,1]$ continuously to $[0,1)$ and is increasing, as may be verified by differentiation, the equation $\phi(g) = g$ has a root in $[0,1)$ and $g_n$ converges to the largest such root $g^*$ as $n \to \infty$. See Figure 4. Suppose that we accept the branching process as a suitable approximation to the hypergraph for the calculation of the core. Then we are led to the following values for the limiting core frequencies:

$$\bar{p}_{d,d'} = \binom{d'}{d} \sigma(g^*)^d (1 - \sigma(g^*))^{d'-d} p_{d'}, \quad k \leqslant d \leqslant d',$$
$$\bar{q}_w = (g^*)^w q_w, \quad w \geqslant 1. \tag{11}$$

We have not justified the interchange of limits which would be required to turn this into a rigorous argument. This seems unlikely to be straightforward. For, by analogy with Theorem 2.2 in [4], in the critical case when $\phi(g) \leqslant g$ in a neighbourhood of $g^*$, we would expect that, asymptotically, the core frequencies would take values corresponding to smaller roots of $\phi(g) = g$ with probability $1/2$. Thus, in this case, when also the crossing condition of Theorem 7.1 fails, the branching process heuristic would lead to an incorrect conclusion. However, for certain random graphs, this sort of approach was made to work in [20].

### 7.3. Statement of result

We return to the framework described in Subsection 7.1. Thus now $\boldsymbol{n}(\boldsymbol{d}) = N\boldsymbol{p}$ and $\boldsymbol{n}(\boldsymbol{w}) = N\boldsymbol{q}$ for our given frequency vectors $\boldsymbol{p}$ and $\boldsymbol{q}$. Define the distributions $\lambda$ and $\sigma$ by the equations (8)[10]. The normalized core frequencies $\bar{P}_{d,d'}$ and $\bar{Q}_w$ were defined at (7) and the limiting core frequencies $\bar{p}_{d,d'}$ and $\bar{q}_w$ were defined at (11). The following result will be proved in Subsection 7.7 using a differential equation approximation to a suitably chosen Markov chain.

**Theorem 7.1.** *Assume that either $g^* = 0$ or the following crossing condition holds:*
$$g^* = \sup\{g \in [0,1) : \phi(g) > g\}.$$
*Then, for all $\nu \in (0,1]$, there is a constant $C < \infty$, depending only on $\boldsymbol{p}, \boldsymbol{q}$ and $\nu$, such that, for all $N \geqslant 1$,*
$$\mathbb{P}\left(\sup_{k \leqslant d \leqslant d'} |\bar{P}_{d,d'} - \bar{p}_{d,d'}| > \nu \text{ or } \sup_{w \geqslant 1} |\bar{Q}_w - \bar{q}_w| > \nu\right) \leqslant Ce^{-N/C}.$$

### 7.4. Splitting property

A uniform random hypergraph $\Gamma \sim U(\boldsymbol{d}, \boldsymbol{w})$ has a useful *splitting property*, which we now describe. Given a partition $V = V' \cup V''$, we can identify a hypergraph $h$ on $V \times E$ with the pair of hypergraphs $h', h''$ on $V' \times E, V'' \times E$ respectively, obtained by intersection. Consider the partition
$$G(\boldsymbol{d}, \boldsymbol{w}) = \cup_{\boldsymbol{w}'+\boldsymbol{w}''=\boldsymbol{w}} G(\boldsymbol{d}', \boldsymbol{w}') \times G(\boldsymbol{d}'', \boldsymbol{w}''),$$
where $\boldsymbol{d}', \boldsymbol{d}''$ are the restrictions of $\boldsymbol{d}$ to $V', V''$ respectively, and where $\boldsymbol{w}', \boldsymbol{w}''$ range over all weight functions on $E$ subject to the given constraint. We deduce that, conditional on $\{\boldsymbol{W}' = \boldsymbol{w}', \boldsymbol{W}'' = \boldsymbol{w}''\}$, the hypergraphs $\Gamma'$ and $\Gamma''$ are independent, with $\Gamma' \sim U(\boldsymbol{d}', \boldsymbol{w}')$ and $\Gamma'' \sim U(\boldsymbol{d}'', \boldsymbol{w}'')$. By symmetry, an analogous splitting property holds in respect of any partition of $E$. In particular, if $v \in V$ and $e \in E$ are chosen independently of $\Gamma$, then $\Gamma \setminus \Gamma(v)$ and $\Gamma \setminus \Gamma(e)$ are also uniformly distributed given their vertex degrees and edge weights.

---

[10] We have chosen for simplicity to consider a sequential limit in which these distributions remain fixed: the interpretation of (8) in the preceding subsection differs by a factor of $N$, top and bottom, which cancels to leave $\lambda$ and $\sigma$ independent of $N$.

### 7.5. Analysis of a core-finding algorithm

Given a hypergraph $\gamma$ on $V \times E$, set $\gamma_0 = \gamma$ and define recursively a sequence of hypergraphs $(\gamma_n)_{n \geqslant 0}$ as follows: for $n \geqslant 0$, given $\gamma_n$, choose if possible, uniformly at random, a vertex $v_{n+1} \in V$ such that $d = \boldsymbol{d}_{\gamma_n}(v_{n+1}) \in \{1, \ldots, k-1\}$ and set

$$\gamma_{n+1} = \gamma_n \setminus (\gamma_n(e_1) \cup \cdots \cup \gamma_n(e_d)),$$

where $\gamma_n(v_{n+1}) = \{(v_{n+1}, e_1), \ldots, (v_{n+1}, e_d)\}$; if there is no such vertex, set $\gamma_{n+1} = \gamma_n$. Thus we remove from $\gamma_n$ all edges containing the chosen vertex $v_{n+1}$. The sequence terminates at the $k$-core $\bar{\gamma}$.

Take $\Gamma \sim U(\boldsymbol{d}, \boldsymbol{w})$ and consider the corresponding sequence $(\Gamma_n)_{n \geqslant 0}$. We continue to write $v_n$ for the random vertices chosen in the algorithm. Set $\boldsymbol{D}_n = \boldsymbol{d}_{\Gamma_n}$ and $\boldsymbol{W}_n = \boldsymbol{w}_{\Gamma_n}$. In the sequel we shall use the symbols $j, k, l, d, d'$ to denote elements of $\mathbb{Z}^+ \times \{V\}$, while $w$ will denote an element of $\mathbb{Z}^+ \times \{E\}$. This is just a formal device which will allow us to refer to two different sets of coordinates by $\xi^d$ and $\xi^w$, and, to lighten the notation, we shall identify both these sets with $\mathbb{Z}^+$ where convenient. For $0 \leqslant d \leqslant d'$ and $w \geqslant 0$, set

$$\xi_n^{d,d'} = |\{v \in V : \boldsymbol{D}_n(v) = d, \boldsymbol{d}(v) = d'\}|, \quad \xi_n^w = |\{e \in E : \boldsymbol{W}_n(e) = w\}|.$$

Set

$$\xi_n = (\xi_n^{d,d'}, \xi_n^w : 0 \leqslant d \leqslant d', w \geqslant 0).$$

Note that the process $(\xi_n)_{n \geqslant 0}$ is adapted to the filtration $(\mathcal{F}_n)_{n \geqslant 0}$ given by

$$\mathcal{F}_n = \sigma(\boldsymbol{D}_r, \boldsymbol{W}_r : r = 0, 1, \ldots, n).$$

**Proposition 7.2.** *For all $n \geqslant 0$, conditional on $\mathcal{F}_n$, we have $\Gamma_n \sim U(\boldsymbol{D}_n, \boldsymbol{W}_n)$.*

*Proof.* The claim is true for $n = 0$ by assumption. Suppose inductively that the claim holds for $n$. The algorithm terminates on the $\mathcal{F}_n$-measurable event

$$\{\boldsymbol{D}_n(v) \in \{0\} \cup \{k, k+1, \ldots\} \text{ for all } v \in V\},$$

so on this event the claim holds also for $n+1$. Suppose then that the algorithm does not terminate at $n$. Conditional on $\mathcal{F}_n$, $v_{n+1}$ and $\Gamma_n$ are independent. Hence, by splitting, $\Gamma_n \setminus \Gamma_n(v_{n+1})$ is uniform given its vertex degrees and edge weights. Then, by a further splitting, we can delete each of the edges $\Gamma_n(e)$ with $(v_{n+1}, e) \in \Gamma_n$, still preserving this uniform property, to obtain $\Gamma_{n+1}$. Hence the claim holds for $n+1$ and the induction proceeds. □

Note that the conditional distribution of $v_{n+1}$ given $\mathcal{F}_n$ depends only on $\boldsymbol{D}_n$ and that $(\boldsymbol{D}_{n+1}, \boldsymbol{W}_{n+1})$ is a function of $\Gamma_{n+1}$, and hence is a function of $(v_{n+1}, \Gamma_n)$. It follows that $(\boldsymbol{D}_n, \boldsymbol{W}_n)_{n \geqslant 0}$ is a Markov chain and hence, by symmetry, $(\xi_n)_{n \geqslant 0}$ is also a Markov chain. It will be convenient to denote the state-space by $S$, to define for $\xi \in S$,

$$\xi^d = \sum_{d'=d}^{L} \xi^{d,d'}, \quad n(\xi) = \sum_{d=1}^{k-1} \xi^d, \quad l(\xi) = \sum_{d=1}^{k-1} d\xi^d, \quad h(\xi) = \sum_{d=k}^{L} d\xi^d,$$

and
$$m(\xi) = \sum_{w=1}^{L} w\xi^w, \quad p(\xi) = \sum_{w=1}^{L} w(w-1)\xi^w,$$

and to set $q(\xi) = m(\xi)n(\xi)/l(\xi)$. Thus, $\xi^d$ is the number of vertices of degree $d$, and $n(\xi)$ is the number of *light* vertices, that is, those of degree less than $k$; $l(\xi)$ is the total degree of the light vertices, and $h(\xi)$ is the total degree of the *heavy* vertices; $m(\xi)$ is the total weight, and $p(\xi)$ is the number of ordered pairs of elements of $\xi$ having the same edge label. Note that, for all $\xi \in S$, $m(\xi) \leqslant Nm$ and $n(\xi) \leqslant l(\xi)$, so $q(\xi) \leqslant Nm$. We obtain a continuous-time Markov chain $(X_t)_{t \geqslant 0}$ by taking $(\xi_n)_{n \geqslant 0}$ as jump chain and making jumps at rate $q(X_t)$. As we saw in Subsection 5.2, in the study of terminal values, we are free to choose a convenient jump rate, which should, in particular ensure that the terminal time remains tight in the limit of interest. Our present choice will have this property. However, it has been chosen also so that the limiting differential equation has a simple form. Define now coordinate functions $x^{d,d'}$ and $x^w$ on $S$, for $k \leqslant d \leqslant d'$ and $w \geqslant 1$, by
$$x^{d,d'}(\xi) = \xi^{d,d'}/N, \quad x^w(\xi) = \xi^w/N.$$

Set
$$\boldsymbol{X}_t = \boldsymbol{x}(X_t) = (x^{d,d'}(X_t), x^w(X_t) : k \leqslant d \leqslant d' \leqslant L, 1 \leqslant w \leqslant L).$$

We consider $\boldsymbol{X}$ as a process in $\mathbb{R}^D$, where $D = \frac{1}{2}(L-k+1)(L-k+2) + L$. We shall use $h(x)$, $m(x)$ and $p(x)$ to denote functions of
$$x = (x^{d,d'}, x^w : k \leqslant d \leqslant d' \leqslant L, 1 \leqslant w \leqslant L) \in \mathbb{R}^D,$$

defined as for $\xi \in S$, but replacing $\xi^{d,d'}$ and $\xi^w$ by $x^{d,d'}$ and $x^w$ respectively. Note that the jumps of $\boldsymbol{X}$ are bounded in supremum norm by $(k-1)(L-1)/N$. Note also that $h(\boldsymbol{X}_t) \leqslant m(\boldsymbol{X}_t)$ for all $t$ and that the algorithm terminates at $T_0 = \inf\{t \geqslant 0 : h(\boldsymbol{X}_t) = m(\boldsymbol{X}_t)\}$. Hence $\boldsymbol{X}_{T_0}$ is the desired vector of core frequencies:
$$\bar{P}_{d,d'} = \boldsymbol{X}_{T_0}^{d,d'}, \quad \bar{Q}_w = \boldsymbol{X}_{T_0}^w. \tag{12}$$

Recall that $m = m(\boldsymbol{p}) = m(\boldsymbol{q})$ is a given constant. We also write $m(x)$ for the function on $\mathbb{R}^D$ just defined. Thus $m = m(\boldsymbol{X}_0)$. Let
$$U_0 = \{x \in \mathbb{R}^D : x^{d,d'}, x^w \in [0,m], m(x) > 0\}$$

and note that $\boldsymbol{x}(\xi) \in U_0$ for all $\xi \in S \setminus \{0\}$. Define a vector field $b$ on $U_0$ by
$$b^{d,d'}(x) = \frac{p(x)}{m(x)}\{(d+1)x^{d+1,d'} - dx^{d,d'}\}, \quad k \leqslant d \leqslant d' \leqslant L,$$

where $x^{d+1,d} = 0$ for $k \leqslant d \leqslant L$, and
$$b^w(x) = -wx^w, \quad 1 \leqslant w \leqslant L.$$

Define, as in Section 4, the drift vector $\beta$ on $S$ by
$$\beta(\xi) = \sum_{\xi' \neq \xi} (\boldsymbol{x}(\xi') - \boldsymbol{x}(\xi))q(\xi,\xi').$$

**Proposition 7.3.** *There is a decreasing function $\psi : \mathbb{N} \to [0,1]$, depending only on $\boldsymbol{p}$ and $\boldsymbol{q}$, with $\psi(N) \to 0$ as $N \to \infty$, such that, for all $\xi \in S$ with $\boldsymbol{x}(\xi) \in U_0$,*
$$\|\beta(\xi) - b(\boldsymbol{x}(\xi))\| \leqslant \psi(m(\xi)).$$

*Proof.* Fix $\xi \in S$ and condition on $\xi_0 = \xi$. Then, for $l = 1, \ldots, k-1$, we have $\boldsymbol{d}(v_1) = l$ with probability $\xi^l/n(\xi)$. Condition further on $v_1 = v$ and $\boldsymbol{d}(v) = l$ and write $\Gamma(v) = \{v\} \times \{e_1, \ldots, e_l\}$. We use the notation of Subsection 7.2 for the local structure. Then we have $\xi_1^w - \xi_0^w = -\sum_{i=1}^l 1_{\{S_i = w-1\}}$, so
$$\beta^w(\xi) = q(\xi)\mathbb{E}(\xi_1^w - \xi_0^w | \xi_0 = \xi)/N = -\frac{q(\xi)}{Nn(\xi)} \sum_{l=1}^{k-1} l\xi^l \sigma_{w-1}(\xi, l),$$
where
$$\sigma_{w-1}(\xi, l) = \mathbb{P}(S_1 = w-1 | \xi_0 = \xi, \boldsymbol{d}(v_1) = l).$$
By the branching process approximation, we can find a function $\psi$, of the required form, such that
$$m|\sigma_{w-1}(\xi, l) - w\xi^w/m(\xi)| \leqslant \psi(m(\xi)), \quad w = 1, \ldots, L.$$
After some straightforward estimation we obtain, for the same function $\psi$, the required estimate
$$|\beta^w(\xi) - b^w(\boldsymbol{x}(\xi))| \leqslant \psi(m(\xi)).$$

We turn to the remaining components. Note that $|\xi_1^{d,d'} - \xi_0^{d,d'}| \leqslant (k-1)(L-1)$. Recall from Subsection 7.2 the event
$$A = \{(v', e_i), (v', e_j) \in \Gamma \text{ implies } v' = v \text{ or } i = j\}.$$
Condition on $S_1, \ldots, S_l$ and on $L_{i,j}$ for $j = 1, \ldots, S_i$. On $A$, by symmetry, we have $\xi_1^{d,d'} - \xi_0^{d,d'} = Z^{d+1,d'} - Z^{d,d'}$, where $Z^{d,d'}$ has binomial distribution with parameters $\sum_{i=1}^l \sum_{j=1}^{S_i} 1_{\{L_{i,j} = d-1\}}$ and $\xi^{d,d'}/\xi^d$. Now,
$$\beta^{d,d'}(\xi) = q(\xi)\mathbb{E}(\xi_1^{d,d'} - \xi_0^{d,d'} | \xi_0 = \xi)/N$$
and
$$\mathbb{E}(Z^{d,d'}|\xi_0 = \xi) = \sum_{l=1}^{k-1}(\xi^l/n(\xi))l\sum_{w=1}^L (w-1)\sigma_{w-1}(\xi,l)\lambda_{d-1}(\xi,l,w-1)(\xi^{d,d'}/\xi^d),$$
where
$$\lambda_{d-1}(\xi, l, w-1) = \mathbb{P}(L_{1,1} = d-1 | \xi_0 = \xi, \boldsymbol{d}(v_1) = l, S_1 = w-1).$$
By the branching process approximation, we can find a function $\psi$, of the required form, such that $m\mathbb{P}(A^c) \leqslant \psi(m(\xi))$ and
$$m|\lambda_d(\xi, l, w-1) - (d+1)\xi^{d+1}/m(\xi)| \leqslant \frac{1}{2}(L-1)((L+4k)\psi(m(\xi)).$$
Then, by some straightforward estimation, for the same function $\psi$,
$$|\beta^{d,d'}(\xi) - b^{d,d'}(\boldsymbol{x}(\xi))| \leqslant \psi(m(\xi)).$$
□

### 7.6. Solving the differential equation

Consider the limiting differential equation $\dot{x}_t = b(x_t)$ in $U_0$, with starting point $x_0 = \boldsymbol{X}_0$ given by

$$x_0^w = q_w, \quad x_0^{d,d} = p_d, \quad x_0^{d,d'} = 0, \quad 1 \leqslant w \leqslant L, \quad k \leqslant d < d' \leqslant L.$$

In components, the equation is written

$$\dot{x}_t^w = -wx_t^w, \qquad\qquad 1 \leqslant w \leqslant L,$$

$$\dot{x}_t^{d,d'} = \frac{p(x_t)}{m(x_t)}\{(d+1)x_t^{d+1,d'} - dx_t^{d,d'}\}, \qquad k \leqslant d \leqslant d' \leqslant L.$$

There is a unique solution $(x_t)_{t \geqslant 0}$ in $U_0$ and, clearly, $x_t^w = e^{-tw}q_w$. Then $m(x_t) = me^{-t}\sigma(e^{-t})$ and $p(x_t) = me^{-2t}\sigma'(e^{-t})$. Hence, if $(\tau_t)_{t \geqslant 0}$ is defined by

$$\dot{\tau}_t = p(x_t)/m(x_t), \quad \tau_0 = 0,$$

then $e^{-\tau} = \sigma(e^{-t})$. A straightforward computation now shows that the remaining components of the solution are given by

$$x_t^{d,d'} = \binom{d'}{d}\sigma(e^{-t})^d(1-\sigma(e^{-t}))^{d'-d}p_{d'}$$

and that $h(x_t) = m\phi(e^{-t})\sigma(e^{-t})$. Note that $(m-h)(x_t) = \sigma(e^{-t})(e^{-t} - \phi(e^{-t}))$, so $g^* = e^{-\zeta_0}$, where $\zeta_0 = \inf\{t \geqslant 0 : m(x_t) = h(x_t)\}$.

### 7.7. Proof of Theorem 7.1

Recall that the core frequencies are found at the termination of the core-finding algorithm, see (12). A suitably chosen vector of frequencies evolves under this algorithm as a Markov chain, which we can approximate using the differential equation whose solution we have just obtained. The accuracy of this approximation is good so long as the hypergraph remains large.

Consider first the case where $g^* = 0$, when we have $m(x_t) > h(x_t)$ for all $t \geqslant 0$. Here the hypergraph may become small as the algorithm approaches termination, so we run close to termination and then use a monotonicity argument. Fix $\nu \in (0,1]$, set $\mu = \nu/3$ and choose $t_0$ such that $m(x_{t_0}) = 2\mu$. Define

$$U = \{x \in U_0 : m(x) > h(x) \vee \mu\}$$

and set

$$\zeta = \inf\{t \geqslant 0 : x_t \notin U\}, \quad T = \inf\{t \geqslant 0 : \boldsymbol{X}_t \notin U\},$$

as in Section 4. Since $m(x_t)$ is decreasing in $t$, we have $x_t \in U$ for all $t \leqslant t_0$. Hence there exists $\varepsilon \in (0, \nu/(3L))$, depending only on $\boldsymbol{p}, \boldsymbol{q}$ and $\nu$ such that

$$\text{for all } \xi \in S \text{ and } t \leqslant t_0, \quad \|\boldsymbol{x}(\xi) - x_t\| \leqslant \varepsilon \implies \boldsymbol{x}(\xi) \in U.$$

It is straightforward to check, by bounding the first derivative, that $b$ is Lipschitz on $U$ with constant $K \leqslant (L-1)L^3 m/\mu$. Set $\delta = \varepsilon e^{-Kt_0}/3$, as in Section 4. We have $\boldsymbol{X}_0 = x_0$, so $\Omega_0 = \Omega$. By Proposition 7.3, and using the fact that $m(X_t)$ does not increase, we have

$$\int_0^{T \wedge t_0} \|\beta(X_t) - b(\boldsymbol{x}(X_t))\| dt \leqslant \psi(m(X_T))t_0 = \psi(Nm(\boldsymbol{X}_T))t_0 \leqslant \psi(N\mu)t_0,$$

so $\Omega_1 = \Omega$ provided $N$ is large enough that $\psi(N\mu)t_0 \leqslant \delta$. The total jump rate $q(\xi)$ is bounded by $Q = Nm$ for all $\xi \in S$. The norm of the largest jump is bounded by $J = (k-1)(L-1)/N$. Take $A = QJ^2 e = (k-1)^2(L-1)^2 me/N$ and note that $\delta J/(At_0) \leqslant \delta/((k-1)(L-1)met_0) \leqslant 1$, so $A \geqslant QJ^2 \exp\{\delta J/(At_0)\}$, and so $\Omega_2 = \Omega$ as in Subsection 4.2, Footnote 4. On the event $\{\sup_{t \leqslant t_0} \|\boldsymbol{X}_t - x_t\| \leqslant \varepsilon\}$, we have $T_0 > t_0$, so

$$\sum_{k \leqslant d \leqslant d'} d\bar{P}_{d,d'} = \sum_{w \geqslant 1} w\bar{Q}_w = m(\boldsymbol{X}_{T_0}) \leqslant m(\boldsymbol{X}_{t_0})$$

$$\leqslant m(x_{t_0}) + |m(x_{t_0}) - m(\boldsymbol{X}_{t_0})| \leqslant 2\mu + L\varepsilon \leqslant \nu.$$

Hence, by Theorem 4.2, we obtain

$$\mathbb{P}\left(\sup_{k \leqslant d \leqslant d'} \bar{P}_{d,d'} > \nu \text{ or } \sup_{w \geqslant 1} \bar{Q}_w > \nu\right)$$

$$\leqslant \mathbb{P}\left(\sup_{t \leqslant t_0} \|\boldsymbol{X}_t - x_t\| > \varepsilon\right) \leqslant 2De^{-\delta^2/(2At_0)} \leqslant Ce^{-N/C},$$

for a constant $C \in [1, \infty)$ depending only on $\boldsymbol{p}, \boldsymbol{q}$ and $\nu$, which is the conclusion of the theorem in the case $g^* = 0$.

We turn to the case where $g^* > 0$ and $g^* = \sup\{g \in [0,1) : \phi(g) > g\}$. Set now $\mu = \frac{1}{2}m(x_{\zeta_0})$ and choose $t_0 > \zeta_0$. Define $U$, $\zeta$ and $T$ as in the preceding paragraph, noting that $\zeta = \zeta_0$. We seek to apply the refinement of Theorem 4.3 described in Footnote 6, and refer to Subsection 4.3 for the definition of $\rho$. By the crossing condition, $\phi(g) > g$ immediately below $g^*$, so $(m-h)(x_t) = \sigma(e^{-t})(e^{-t} - \phi(e^{-t})) < 0$ immediately after $\zeta = -\log g^*$. We have

$$|(m-h)(x) - (m-h)(x')| \leqslant C\|x - x'\|, \quad |m(x) - m(x')| \leqslant C\|x - x'\|,$$

for a constant $C < \infty$ depending only on $L$. So, given $\nu > 0$, we can choose $\varepsilon > 0$, depending only on $\boldsymbol{p}, \boldsymbol{q}$ and $\nu$, such that $\varepsilon + \rho(\varepsilon) \leqslant \nu$ and $C(\varepsilon + \rho(\varepsilon)) < \frac{1}{2}m(x_\zeta)$. Note that $\|\boldsymbol{X}_T - x_\zeta\| \leqslant \varepsilon + \rho(\varepsilon)$ implies that $m(X_T) > \frac{1}{2}m(x_\zeta)$ and hence that $T = T_0$. Define $\delta$ and $A$ as in the preceding paragraph. Then, by a similar argument, provided $N$ is sufficiently large that $\psi(N\mu)t_0 \leqslant \delta$, we have $\Omega_0 = \Omega_1 = \Omega_2 = \Omega$. Hence, by Theorem 4.3,

$$\mathbb{P}\left(\sup_{k \leqslant d \leqslant d'} |\bar{P}_{d,d'} - \bar{p}_{d,d'}| > \nu \text{ or } \sup_{w \geqslant 1} |\bar{Q}_w - \bar{q}_w| > \nu\right)$$

$$= \mathbb{P}\left(\|\boldsymbol{X}_T - x_\zeta\| > \varepsilon + \rho(\varepsilon)\right) \leqslant 2De^{-\delta^2/(2At_0)} \leqslant Ce^{-N/C},$$

for a constant $C \in [1, \infty)$ depending only on $\boldsymbol{p}, \boldsymbol{q}$ and $\nu$, as required.

## 8. Appendix: Identification of martingales for a Markov chain

We discuss in this appendix the identification of martingales associated with a continuous-time Markov chain $X = (X_t)_{t \geq 0}$ with finite jump rates. In keeping with the rest of the paper, we assume that $X$ has a countable state-space, here denoted $E$, and write $Q = (q(x,y) : x, y \in E)$ for the associated generator matrix. An extension to the case of a general measurable state-space is possible and requires only cosmetic changes. A convenient and elementary construction of such a process $X$ may be given in terms of its *jump chain* $(Y_n)_{n \geq 0}$ and *holding times* $(S_n)_{n \geq 1}$. We shall deduce, directly from this construction, a method to identify the martingales associated with $X$, which proceeds by expressing them in terms of a certain integer-valued random measure $\mu$. There is a close analogy between this method and the common use of Itô's formula in the case of diffusion processes. The method is well known to specialists but we believe there is value in this direct derivation from the elementary construction. Our arguments in this section involve more measure theory than the rest of the paper; we do not however need the theory of Markov semigroups.

### 8.1. The jump-chain and holding-time construction

The jump chain is a sequence $(Y_n)_{n \geq 0}$ of random variables in $E$, and the holding times $(S_n)_{n \geq 1}$ are non-negative random variables which may sometimes take the value $\infty$. We specify the distributions of these random variables in terms of the *jump matrix* $\Pi = (\pi(x,y) : x, y \in E)$ and the *jump rates* $(q(x) : x \in E)$, given by

$$\pi(x,y) = \begin{cases} q(x,y)/q(x), & y \neq x \text{ and } q(x) \neq 0, \\ 0, & y \neq x \text{ and } q(x) = 0, \end{cases} \qquad \pi(x,x) = \begin{cases} 0, & q(x) \neq 0, \\ 1, & q(x) = 0, \end{cases}$$

and

$$q(x) = -q(x,x) = \sum_{y \neq x} q(x,y).$$

Take $Y = (Y_n)_{n \geq 0}$ to be a discrete-time Markov chain with transition matrix $\Pi$. Thus, for all $n \geq 0$, and all $x_0, x_1, \ldots, x_n \in E$,

$$\mathbb{P}(Y_0 = x_0, Y_1 = x_1, \ldots, Y_n = x_n) = \lambda(x_0)\pi(x_0, x_1) \ldots \pi(x_{n-1}, x_n),$$

where $\lambda(x) = \mathbb{P}(Y_0 = x)$. Take $(T_n)_{n \geq 1}$ to be a sequence of independent exponential random variable of parameter 1, independent of $Y$. Set

$$S_n = T_n/q(Y_{n-1}), \quad J_0 = 0, \quad J_n = S_1 + \cdots + S_n, \quad \zeta = \sum_{n=1}^{\infty} S_n,$$

and construct $X$ by

$$X_t = \begin{cases} Y_n, & J_n \leq t < J_{n+1}, \\ \partial, & t \geq \zeta, \end{cases}$$

where $\partial$ is some *cemetery state*, which we adjoin to $E$. We are now using $\zeta$ for the explosion time of the Markov chain $X$, at variance with the rest of the paper. For $t \geqslant 0$, define

$$J_1(t) = \inf\{s \geqslant t : X_s \neq X_t\}, \quad Y_1(t) = X_{J_1(t)}.$$

These are, respectively, the time and destination of the first jump of $X$ starting from time $t$. Consider the *natural filtration* $(\mathcal{F}_t^X)_{t \geqslant 0}$, given by $\mathcal{F}_t^X = \sigma(X_s : s \leqslant t)$. Write $\mathcal{E}$ for the set of subsets of $E$ and set $q(\partial) = 0$ and $\pi(x, B) = \sum_{y \in B} \pi(x, y)$ for $B \in \mathcal{E}$.

**Proposition 8.1.** *For all $s, t \geqslant 0$ and all $B \in \mathcal{E}$, we have, almost surely,*

$$\mathbb{P}(J_1(t) > t + s, Y_1(t) \in B | \mathcal{F}_t^X) = \pi(X_t, B) e^{-q(X_t)s}.$$

Before proving the proposition, we need a lemma, which expresses in precise terms that, if $X$ has made exactly $n$ jumps by time $t$, then all we know at that time are the states $Y_0, \ldots, Y_n$, the times $J_1, \ldots, J_n$ and the fact that the next jump happens later.

**Lemma 8.2.** *Define $\mathcal{G}_n = \sigma(Y_m, J_m : m \leqslant n)$. For all $A \in \mathcal{F}_t^X$ and all $n \geqslant 0$, there exists $\tilde{A}_n \in \mathcal{G}_n$ such that*

$$A \cap \{J_n \leqslant t < J_{n+1}\} = \tilde{A}_n \cap \{t < J_{n+1}\}.$$

*Proof.* Denote by $\mathcal{A}_t$ the set of all sets $A \in \mathcal{F}_t^X$ for which the desired property holds. Then $\mathcal{A}_t$ is a $\sigma$-algebra. For any $s \leqslant t$, we can write

$$\{X_s \in B\} \cap \{J_n \leqslant t < J_{n+1}\} = \tilde{A}_n \cap \{t < J_{n+1}\},$$

where

$$\tilde{A}_n = \cup_{m=0}^{n-1}\{Y_m \in B, J_m \leqslant s < J_{m+1}, J_n \leqslant t\} \cup \{Y_n \in B, J_n \leqslant s\},$$

so $\{X_s \in B\} \in \mathcal{A}_t$. Hence $\mathcal{A}_t = \mathcal{F}_t^X$. □

*Proof of Proposition 8.1.* The argument relies on the memoryless property of the exponential distribution, in the following conditional form: for $s, t \geqslant 0$ and $n \geqslant 0$, almost surely, on $\{J_n \leqslant t\}$,

$$\mathbb{P}(J_{n+1} > t + s | \mathcal{G}_n) = \mathbb{P}(T_{n+1} > q(Y_n)(s + t - J_n) | \mathcal{G}_n)$$
$$= e^{-q(Y_n)(s+t-J_n)} = e^{-q(Y_n)s} \mathbb{P}(J_{n+1} > t | \mathcal{G}_n).$$

Then for $B \in \mathcal{E}$ and $A \in \mathcal{F}_t^X$, we have

$$\mathbb{P}(J_1(t) > t + s, Y_1(t) \in B, A, J_n \leqslant t < J_{n+1}) = \mathbb{P}(J_{n+1} > t + s, Y_{n+1} \in B, \tilde{A}_n)$$
$$= \mathbb{E}(\pi(Y_n, B) e^{-q(Y_n)s} 1_{\tilde{A}_n \cap \{J_{n+1} > t\}}) = \mathbb{E}(\pi(X_t, B) e^{-q(X_t)s} 1_{A \cap \{J_n \leqslant t < J_{n+1}\}})$$

and

$$\mathbb{P}(J_1(t) > t + s, Y_1(t) \in B, A, t \geqslant \zeta) = \delta_\partial(B) = \mathbb{E}(\pi(X_t, B) e^{-q(X_t)s} 1_{A \cap \{t \geqslant \zeta\}})$$

On summing all the above equations we obtain

$$\mathbb{P}(J_1(t) > t + s, Y_1(t) \in B, A) = \mathbb{E}(\pi(X_t, B) e^{-q(X_t)s} 1_A),$$

as required. □

### 8.2. Markov chains in a given filtration

For many purposes, the construction of a process $X$ which we have just given serves as a good definition of a continuous-time Markov chain with generator $Q$. However, from now on, we adopt a more general definition, which has the merit of expressing a proper relationship between $X$ and a general given filtration $(\mathcal{F}_t)_{t \geqslant 0}$. Assume that $X$ is *constant on the right*, that is to say, for all $t \geqslant 0$, there exists $\varepsilon > 0$ such that $X_s = X_t$ whenever $t \leqslant s < t + \varepsilon$. Set $J_0 = 0$ and define for $n \geqslant 0$,

$$Y_n = X_{J_n}, \quad J_{n+1} = \inf\{t \geqslant J_n : X_t \neq X_{J_n}\}. \tag{13}$$

For $t \geqslant 0$, define $J_1(t)$ and $Y_1(t)$ as above. Assume that $X$ is *minimal*, so that $X_t = \partial$ for all $t \geqslant \zeta$, where $\zeta = \sup_n J_n$. Assume finally that $X$ is adapted to $(\mathcal{F}_t)_{t \geqslant 0}$. Then we say that $X$ is a *continuous-time $(\mathcal{F}_t)_{t \geqslant 0}$-Markov chain with generator $Q$* if, for all $s, t \geqslant 0$, and all $B \in \mathcal{E}$, we have, almost surely,

$$\mathbb{P}(J_1(t) > t + s, Y_1(t) \in B | \mathcal{F}_t) = \pi(X_t, B) e^{-q(X_t)s}.$$

The process constructed above from jump chain and holding times is constant on the right and minimal and we do recover the jump chain and holding times using (13); moreover by Proposition 8.1, such a process is then a continuous-time Markov chain in its natural filtration. The defining property of a continuous-time Markov chain extends to stopping times.

**Proposition 8.3.** *Let $X$ be an $(\mathcal{F})_{t \geqslant 0}$-Markov chain with generator $Q$ and let $T$ be a stopping time. Then, for all $s \geqslant 0$ and $B \in \mathcal{E}$, on $\{T < \infty\}$, almost surely,*

$$\mathbb{P}(J_1(T) > T + s, Y_1(T) \in B | \mathcal{F}_T) = \pi(X_T, B) e^{-q(X_T)s}.$$

*Proof.* Consider the stopping times $T_m = 2^{-m} \lceil 2^m T \rceil$. Note that $T_m \downarrow T$ as $m \to \infty$ so, since $X$ is constant on the right, $X_{T_m} = X_T$, $J_1(T_m) = J_1(T)$ and $Y_1(T_m) = Y_1(T)$ eventually as $m \to \infty$, almost surely. Suppose $A \in \mathcal{F}_T$ with $A \subseteq \{T < \infty\}$. Then for all $k \in \mathbb{Z}^+$, $A \cap \{T_m = k2^{-m}\} \in \mathcal{F}_{k2^{-m}}$, so, almost surely,

$$\mathbb{P}(J_1(T_m) > T_m + s, Y_1(T_m) \in B, A, T_m = k2^{-m})$$
$$= \mathbb{E}(\pi(X_{T_m}, B) e^{-q(X_{T_m})s} 1_{A \cap \{T_m = k2^{-m}\}}),$$

and, summing over $k$,

$$\mathbb{P}(J_1(T_m) > T_m + s, Y_1(T_m) \in B, A) = \mathbb{E}(\pi(X_{T_m}, B) e^{-q(X_{T_m})s} 1_A).$$

Letting $m \to \infty$, we can replace, by bounded convergence, $T_m$ by $T$, thus proving the proposition. □

### 8.3. The jump measure and its compensator

The *jump measure* $\mu$ of $X$ and its *compensator* $\nu$ are random measures on $(0, \infty) \times E$, given by

$$\mu = \sum_{t:X_t \neq X_{t-}} \delta_{(t,X_t)} = \sum_{n=1}^{\infty} \delta_{(J_n, Y_n)}$$

and

$$\nu(dt, B) = q(X_{t-}, B)dt = q(X_{t-})\pi(X_{t-}, B)dt, \quad B \in \mathcal{E}.$$

Recall that the previsible $\sigma$-algebra $\mathcal{P}$ on $\Omega \times (0, \infty)$ is the $\sigma$-algebra generated by all left-continuous adapted processes. Extend this notion in calling a function defined on $\Omega \times (0, \infty) \times E$ *previsible* if it is $\mathcal{P} \otimes \mathcal{E}$-measurable.

**Theorem 8.4.** *Let $H$ be previsible and assume that, for all $t \geqslant 0$,*

$$\mathbb{E} \int_0^t \int_E |H(s, y)| \nu(ds, dy) < \infty.$$

*Then the following process is a well-defined martingale*

$$M_t = \int_{(0,t] \times E} H(s, y)(\mu - \nu)(ds, dy).$$

Define measures $\bar{\mu}$ and $\bar{\nu}$ on $\mathcal{P} \otimes \mathcal{E}$ by

$$\bar{\mu}(D) = \mathbb{E}(\mu(D)), \quad \bar{\nu}(D) = \mathbb{E}(\nu(D)), \quad D \in \mathcal{P} \otimes \mathcal{E}.$$

We shall show that $\bar{\mu} = \bar{\nu}$. Once this is done, the proof of Theorem 8.4 will be straightforward. For $n \geqslant 0$, define measures $\bar{\mu}_n$ and $\bar{\nu}_n$ on $\mathcal{P} \otimes \mathcal{E}$ by

$$\bar{\mu}_n(D) = \bar{\mu}(D \cap (J_n, J_{n+1}]), \quad \bar{\nu}_n(D) = \bar{\nu}(D \cap (J_n, J_{n+1}]), \quad D \in \mathcal{P} \otimes \mathcal{E},$$

where $D \cap (J_n, J_{n+1}] = \{(\omega, t, y) \in D : J_n(\omega) < t \leqslant J_{n+1}(\omega)\}$. Then, since $q(\partial) = 0$,

$$\bar{\mu} = \sum_{n=0}^{\infty} \bar{\mu}_n, \quad \bar{\nu} = \sum_{n=0}^{\infty} \bar{\nu}_n,$$

so it will suffice to show the following lemma.

**Lemma 8.5.** *For all $n \geqslant 0$, we have $\bar{\mu}_n = \bar{\nu}_n$.*

*Proof.* The proof rests on the following basic identity for an exponential random variable $V$ of parameter $q$:

$$\mathbb{P}(V \leqslant s) = q\mathbb{E}(V \wedge s), \quad s \geqslant 0.$$

Let $T$ be a stopping time and let $S$ be a non-negative $\mathcal{F}_T$-measurable random variable. Set $U = (T + S) \wedge J_1(T)$. By Proposition 8.3, we know that, conditional on $\mathcal{F}_T$, $J_1(T)$ and $Y_1(T)$ are independent, $J_1(T)-T$ has exponential distribution of parameter $q(X_T)$ and $Y_1(T)$ has distribution $\pi(X_T,.)$. From the basic identity, we obtain

$$\mathbb{P}(J_1(T) - T \leqslant S|\mathcal{F}_T) = q(X_T)\mathbb{E}((J_1(T) - T) \wedge S|\mathcal{F}_T) = q(X_T)\mathbb{E}(U - T|\mathcal{F}_T).$$

Fix $n \geqslant 0$, $t \leqslant u$, $A \in \mathcal{F}_t$, $B \in \mathcal{E}$ and set $D = A \times (t,u] \times B$. The set of such sets $D$ forms a $\pi$-system, which generates the $\sigma$-algebra $\mathcal{P} \otimes \mathcal{E}$. We shall show that $\bar{\mu}_n(D) = \bar{\nu}_n(D) \leqslant 1$. By taking $A = \Omega$, $B = E$, $t = 0$ and letting $u \to \infty$, this shows also that $\bar{\mu}_n$ and $\bar{\nu}_n$ have the same finite total mass. The lemma will then follow by uniqueness of extension.

Take $T = J_n \wedge t \vee J_{n+1}$ and $S = (u-T)^+ 1_{\{T < J_{n+1}\}}$. Then

$$U = (T + S) \wedge J_1(T) = J_n \wedge u \vee J_{n+1}.$$

So

$$\begin{aligned}
\bar{\mu}_n(D) &= \bar{\mu}(D \cap (J_n, J_{n+1}]) = \mathbb{P}(J_1(T) \leqslant T + S, Y_1(T) \in B, A) \\
&= \mathbb{E}(1_A \mathbb{P}(J_1(T) - T \leqslant S, Y_1(T) \in B|\mathcal{F}_T)) \\
&= \mathbb{E}(1_A q(X_T)\pi(X_T, B)\mathbb{E}(U - T|\mathcal{F}_T)) \\
&= \mathbb{E}(1_A q(X_T, B)(U - T)) = \mathbb{E}\left(1_A \int_T^U q(X_s, B)ds\right) \\
&= \bar{\nu}(D \cap (J_n, J_{n+1}]) = \bar{\nu}_n(D),
\end{aligned}$$

as required. $\square$

*Proof of Theorem 8.4.* For a non-negative previsible function $H$, for $s \leqslant t$ and $A \in \mathcal{F}_s$, by Fubini's theorem,

$$\begin{aligned}
\mathbb{E}\left(1_A \int_{(s,t] \times E} H(r,y)\mu(dr,dy)\right) &= \int_{A \times (s,t] \times E} H d\bar{\mu} \\
&= \int_{A \times (s,t] \times E} H d\bar{\nu} = \mathbb{E}\left(1_A \int_{(s,t] \times E} H(r,y)\nu(dr,dy)\right).
\end{aligned}$$

So, taking $A = \Omega$ and $s = 0$, if

$$\mathbb{E}\int_0^t \int_E H(r,y)\nu(dr,dy) < \infty$$

for all $t \geqslant 0$, then $M_t$ is well-defined and integrable and, now with general $s \leqslant t$ and $A \in \mathcal{F}_s$,

$$\mathbb{E}((M_t - M_s)1_A) = 0.$$

The result extends to general previsible functions $H$ by taking differences. $\square$

### 8.4. Martingale estimates

Theorem 8.4 makes it possible to identify martingales associated with $X$ in a manner analogous to Itô's formula. We illustrate this by deriving three martingales $M, N$ and $Z$ associated with a given function $f : E \to \mathbb{R}$. The processes $M$ and $N$ depend, respectively, linearly and quadratically on $f$, whereas $Z$ is an exponential martingale. We use $N$ and $Z$ to obtain quadratic and exponential martingale inequalities for $M$, which are used in the main part of the paper. We emphasise that $f$ can be any function. In the main part of the paper, we work with the martingales associated with several choices of $f$ at once. In this subsection we do not burden the notation by registering further the dependence of everything on $f$. The discussion that follows has a computational aspect and an analytic aspect. The reader may wish to check the basic computations before considering in detail the analytic part. We note for orientation that, in the simple case where the maximum jump rate is bounded and where $f$ also is bounded, then there is no explosion and $M, N$ and $Z$, as defined below, are all martingales, without any need for reduction by stopping times. For simplicity, we make an assumption in this subsection that $X$ *does not explode*. A reduction to this case is always possible by an adapted random time-change – this can allow the identification of martingales in the explosive case by applying the results given below and then inverting the time-change. We omit further details.

For all $t \in [0, \infty)$, we have $J_n \leqslant t < J_{n+1}$ for some $n \geqslant 0$. Then

$$f(X_t) = f(Y_n) = f(Y_0) + \sum_{m=0}^{n-1} \{f(Y_{m+1}) - f(Y_m)\}$$
$$= f(X_0) + \int_{(0,t] \times E} \{f(y) - f(X_{s-})\} \mu(ds, dy).$$

Define

$$\tau(x) = \sum_{y \neq x} |f(y) - f(x)| q(x, y)$$

and set $\zeta_1 = \inf\{t \geqslant 0 : \tau(X_t) = \infty\}$. Define when $\tau(x) < \infty$

$$\beta(x) = \sum_{y \neq x} \{f(y) - f(x)\} q(x, y).$$

Then, for $t \in [0, \infty)$ with $t \leqslant \zeta_1$,

$$\int_{(0,t] \times E} \{f(y) - f(X_{s-})\} \nu(ds, dy)$$
$$= \int_0^t \int_E \{f(y) - f(X_{s-})\} q(X_{s-}, dy) ds = \int_0^t \beta(X_s) ds,$$

so

$$f(X_t) = f(X_0) + M_t + \int_0^t \beta(X_s) ds,$$

where
$$M_t = \int_0^t \int_E \{f(y) - f(X_{s-})\}(\mu - \nu)(ds, dy).$$

Define, as usual, for stopping times $T$, the stopped process $M_t^T = M_{T \wedge t}$.

**Proposition 8.6.** *For all stopping times $T \leqslant \zeta_1$, we have*
$$\mathbb{E}\left(\sup_{t \leqslant T} |M_t|\right) \leqslant 2\mathbb{E} \int_0^T \tau(X_t) dt,$$

*and, if the right hand side is finite, then $M^T$ is a martingale. Moreover $M^{\zeta_1}$ is always a local martingale.*

*Proof.* Let $T \leqslant \zeta_1$ be a stopping time, with $\mathbb{E} \int_0^T \tau(X_t) dt < \infty$. Consider the previsible process
$$H_1(t, y) = \{f(y) - f(X_{t-})\} 1_{\{t \leqslant T\}}.$$

Then
$$M_t^T = \int_{(0,t] \times E} H_1(s, y)(\mu - \nu)(ds, dy),$$
$$\sup_{t \leqslant T} |M_t| \leqslant \int_{(0,\infty) \times E} |H_1(t, y)|(\mu + \nu)(dt, dy)$$

and
$$\int_{(0,\infty) \times E} |H_1(t, y)| \nu(dt, dy) = \int_0^T \tau(X_t) dt.$$

The first sentence of the statement now follows easily from Theorem 8.4. For the second, it suffices to note that, for the stopping times $T_n = \inf\{t \geqslant 0 : \tau(X_t) > n\} \wedge n$, we have $T_n \uparrow \zeta_1$ as $n \to \infty$ and $\int_0^{T_n} \tau(X_t) dt \leqslant n^2$, for all $n$, almost surely. □

We turn now to $L^2$ estimates, in the process identifying the martingale decomposition of $M^2$. Note first the following identity: for $t \in [0, \infty)$ with $t \leqslant \zeta_1$,
$$M_t^2 = 2 \int_{(0,t] \times E} M_{s-}\{f(y) - f(X_{s-})\}(\mu - \nu)(ds, dy)$$
$$+ \int_{(0,t] \times E} \{f(y) - f(X_{s-})\}^2 \mu(ds, dy). \tag{14}$$

This may be established by verifying that the jumps of left and right hand sides agree, and that their derivatives agree between jump times. Define
$$\alpha(x) = \sum_{x' \neq x} \{f(x') - f(x)\}^2 q(x, x'),$$

and set $\zeta_2 = \inf\{t \geqslant 0 : \alpha(X_t) = \infty\}$. By Cauchy–Schwarz, we have $\tau(x)^2 \leqslant \alpha(x)q(x)$ for all $x$, so $\zeta_2 \leqslant \zeta_1$. For $t \in [0, \infty)$ with $t \leqslant \zeta_2$ we can write,

$$M_t^2 = N_t + \int_0^t \alpha(X_s)ds, \tag{15}$$

where

$$N_t = \int_{(0,t] \times E} H(s,y)(\mu - \nu)(ds, dy),$$

and

$$H(s,y) = 2M_{s-}\{f(y) - f(X_{s-})\} + \{f(y) - f(X_{s-})\}^2.$$

**Proposition 8.7.** *For all stopping times $T \leqslant \zeta_1$, we have*

$$\mathbb{E}\left(\sup_{t \leqslant T} |M_t|^2\right) \leqslant 4\mathbb{E}\int_0^T \alpha(X_t)dt.$$

*Moreover, $N^{\zeta_2}$ is a local martingale and, for all stopping times $T \leqslant \zeta_2$ with $\mathbb{E}\int_0^T \alpha(X_t)dt < \infty$, both $M^T$ and $N^T$ are martingales, and*

$$\mathbb{E}\left(\sup_{t \leqslant T} |N_t|\right) \leqslant 5\mathbb{E}\int_0^T \alpha(X_t)dt.$$

*Proof.* Let $T \leqslant \zeta_1$ be a stopping time, with $\mathbb{E}\int_0^T \alpha(X_t)dt < \infty$, then $T \leqslant \zeta_2$. Consider the previsible process

$$H_2(t,y) = H(t,y)1_{\{t \leqslant T \wedge T_n\}},$$

where $T_n$ is defined in the preceding proof. Then,

$$N_{T \wedge T_n} = \int_{(0,\infty) \times E} H_2(t,y)(\mu - \nu)(dt, dy)$$

and

$$\mathbb{E}\int_{(0,\infty) \times E} |H_2(s,y)|\nu(ds, dy)$$
$$\leqslant \mathbb{E}\int_0^{T \wedge T_n} \{2|M_t|\tau(X_t) + \alpha(X_t)\}dt \leqslant 4n^4 + \mathbb{E}\int_0^T \alpha(X_t)dt < \infty.$$

Hence, by Theorem 8.4, the process $N^{T \wedge T_n}$ is a martingale. Replace $t$ by $T \wedge T_n$ in (15) and take expectations to obtain

$$\mathbb{E}(|M_{T \wedge T_n}|^2) = \mathbb{E}\int_0^{T \wedge T_n} \alpha(X_t)dt.$$

Apply Doob's $L^2$-inequality to the martingale $M^{T \wedge T_n}$ to obtain

$$\mathbb{E}\left(\sup_{t \leqslant T \wedge T_n} |M_t|^2\right) \leqslant 4\mathbb{E}\int_0^{T \wedge T_n} \alpha(X_t)dt.$$

Go back to (15) to deduce

$$\mathbb{E}\left(\sup_{t \leqslant T \wedge T_n} |N_t|\right) \leqslant 5\mathbb{E}\int_0^{T \wedge T_n} \alpha(X_t)dt.$$

On letting $n \to \infty$, we have $T_n \uparrow \zeta_1$, so we obtain the claimed estimates by monotone convergence, which then imply that $M^T$ and $N^T$ are martingales. Then we can let $T$ run through the sequence $\tilde{T}_n = \inf\{t \geqslant 0 : \alpha(X_t) > n\} \uparrow \zeta_2$ to see that $N^{\zeta_2}$ is a local martingale. □

Finally, we discuss an exponential martingale and estimate. Define for $x \in E$

$$\phi(x) = \sum_{y \neq x} h(f(y) - f(x))q(x, y),$$

where $h(a) = e^a - 1 - a \geqslant 0$, and set $\zeta^* = \inf\{t \geqslant 0 : \phi(X_t) = \infty\}$. Since $e^a - a \geqslant |a|$ for all $a \in \mathbb{R}$, we have $\tau(x) \leqslant \phi(x) + q(x)$, so $\zeta^* \leqslant \zeta_1$. Define for $t \in [0, \infty)$ with $t \leqslant \zeta^*$

$$Z_t = \exp\left\{M_t - \int_0^t \phi(X_s)ds\right\}.$$

Then

$$Z_t = Z_0 + \int_{(0,t] \times E} H^*(s, y)(\mu - \nu)(ds, dy),$$

where

$$H^*(s, y) = Z_{s-}\{e^{f(y) - f(X_{s-})} - 1\}.$$

This identity may be verified in the same way as (14). Consider for $n \geqslant 0$ the stopping time

$$U_n = \inf\{t \geqslant 0 : \phi(X_t) + \tau(X_t) > n\}$$

and note that $U_n \uparrow \zeta^*$ as $n \to \infty$.

**Proposition 8.8.** *For all stopping times* $T \leqslant \zeta_1$,

$$\mathbb{E}\left(\exp\left\{M_T - \int_0^T \phi(X_t)dt\right\}\right) \leqslant 1,$$

*and, for all* $A, B \in [0, \infty)$,

$$\mathbb{P}\left(\sup_{t \leqslant T} M_t > B \text{ and } \int_0^T \phi(X_t)dt \leqslant A\right) \leqslant e^{A-B}.$$

*Moreover,* $Z^{\zeta^*}$ *is a local martingale and a supermartingale, and* $Z^{U_n}$ *is a martingale for all* $n$.

*Proof.* Consider for $n \geqslant 0$ the stopping time

$$V_n = \inf\{t \geqslant 0 : \phi(X_t) + \tau(X_t) > n \text{ or } Z_t > n\}$$

and note that $V_n \uparrow \zeta^*$ as $n \to \infty$. For all $n \geqslant 0$ and $t \geqslant 0$, we have

$$\mathbb{E}\int_{(0,V_n \wedge t] \times E} |H^*(s,y)|\nu(ds,dy) \leqslant \mathbb{E}\int_0^{V_n \wedge t} |Z_s|\{\phi(X_s) + \tau(X_s)\}ds \leqslant n^2 t < \infty.$$

Hence, by Theorem 8.4, for all $n$, the stopped process $Z^{V_n}$ is a martingale. So $Z^{\zeta^*}$ is a local martingale, and hence is a supermartingale by the usual Fatou argument. In particular, for all $t \geqslant 0$, we have $\mathbb{E}(Z_{t \wedge \zeta^*}) \leqslant 1$, so, for all $n \geqslant 0$,

$$\mathbb{E}\int_{(0,U_n \wedge t] \times E} |H^*(s,y)|\nu(ds,dy) \leqslant \mathbb{E}\int_0^{U_n \wedge t} |Z_s|\{\phi(X_s) + \tau(X_s)\}ds \leqslant nt < \infty,$$

and so $Z^{U_n}$ is a martingale. If $T \leqslant \zeta_1$ is a stopping time, then $\mathbb{E}(Z_T) \leqslant \mathbb{E}(Z_{T \wedge \zeta^*})$ and, by optional stopping, $\mathbb{E}(Z_{T \wedge \zeta^*}) \leqslant \mathbb{E}(Z_0) = 1$, so $\mathbb{E}(Z_T) \leqslant 1$, as required. Finally, we can apply this estimate to $T_B = \inf\{t \geqslant 0 : M_t > B\} \wedge T$, noting that $Z_{T_B} \geqslant e^{B-A}$ on the set

$$\Omega_0 = \left\{\sup_{t \leqslant T} M_t > B \text{ and } \int_0^T \phi(X_t)dt \leqslant A\right\},$$

to obtain $e^{B-A}\mathbb{P}(\Omega_0) \leqslant \mathbb{E}(Z_{T_B}) \leqslant 1$, as required. $\square$